\documentclass{amsart}
\usepackage{pstricks,pstcol,pst-plot,pst-node,pst-tree,amssymb}
\usepackage{pstricks-add}

\numberwithin{equation}{section}
\input{tcilatex}

\begin{document}
\title[Regularity of $C^{1}$ smooth surfaces]{Regularity of $C^{1}$ smooth
surfaces with prescribed $p$-mean curvature in the Heisenberg group}
\author{Jih-Hsin Cheng}
\address[Jih-Hsin Cheng and Jenn-Fang Hwang]{ Institute of Mathematics,
Academia Sinica, Taipei, Taiwan, R.O.C.}
\email[Jih-Hsin Cheng]{ cheng@math.sinica.edu.tw}
\urladdr{}
\thanks{}
\author{Jenn-Fang Hwang}
\email[Jenn-Fang Hwang]{ majfh@math.sinica.edu.tw}
\urladdr{}
\thanks{}
\author{Paul Yang}
\address[Paul Yang]{ Department of Mathematics, Princeton University,
Princeton, NJ 08544, U.S.A. }
\email{yang@Math.Princeton.EDU}
\urladdr{}
\thanks{}
\thanks{}

\begin{abstract}
We consider a $C^{1}$ smooth surface with prescribed $p$(or $H$)-mean
curvature in the $3$-dimensional Heisenberg group. Assuming only the
prescribed $p$-mean curvature $H\in C^{0},$ we show that any characteristic
curve is $C^{2}$ smooth and its (line) curvature equals $-H$ in the
nonsingular domain$.$ By introducing characteristic coordinates and invoking
the jump formulas along characteristic curves, we can prove that the
Legendrian (or horizontal) normal gains one more derivative. Therefore the
seed curves are $C^{2}$ smooth. We also obtain the uniqueness of
characteristic and seed curves passing through a common point under some
mild conditions, respectively. These results can be applied to more general
situations.
\end{abstract}

\maketitle

\bigskip



\section{\textbf{Introduction and statement of the results}}

The $p$-minimal (or $X$-minimal, $\boldsymbol{H}$-minimal) surfaces have
been studied extensively in the framework of geometric measure theory (e.g., 
\cite{GN96}, \cite{FSS01}, \cite{Pau01}). Motivated by the isoperimetric
problem in the Heisenberg group, one also studied nonzero constant $p$-mean
curvature surfaces and the regularity problem (e.g., \cite{Pan82}, \cite{CDG}%
, \cite{LM}, \cite{LR}, \cite{RR}, \cite{MR}, \cite{Pauls}). Starting from
the work \cite{CHMY04} (see also \cite{CH04}), we studied the subject from
the viewpoint of partial differential equations and that of differential
geometry. In fact, the notion of $p$-mean curvature ("$p$-" stands for
"pseudohermitian") can be defined for (hyper) surfaces in a pseudohermitian
manifold. The Heisenberg group as a (flat) pseudohermitian manifold is the
simplest model example, and represents a blow-up limit of general
pseudohermitian manifolds. In this paper we will deal with the regularity
problem, in particular, for a $C^{1}$ smooth surface with prescribed $p$%
-mean curvature in the 3-dimensional Heisenberg group $\boldsymbol{H}^{%
\boldsymbol{1}}$. Since our results hold in quite general situations, we
will just start with the general formulation.

Let $\Omega $ be a domain in $R^{m}$ and $u$ $:$ $\Omega $ $\rightarrow $ $R$
be a $W^{1,1}$ function. Let $\vec{F}$ be an arbitrary (say, $L^{1})$ vector
field on $\Omega ,$ and $H$ $\in $ $L^{\infty }(\Omega ).$ In \cite{CHY} we
consider the following energy functional:

\begin{equation*}
\mathcal{F}(u)\equiv \int_{\Omega }|\nabla u+\vec{F}|+Hu
\end{equation*}

\noindent (we omit the Euclidean volume element). When $\vec{F}$ $\equiv $ $%
0,$ this is the energy functional of least gradient. When $m$ $=$ $2$, $\vec{%
F}$ $\equiv $ $(-y,x),$ and $H$ $=$ $0,$ this is the $p($or $\boldsymbol{H)}$%
-area of the graph defined by $u$ in $\boldsymbol{H}^{\boldsymbol{1}}$. Let $%
\varphi $ $\in $ $W^{1,1}(\Omega )$ and $u_{\varepsilon }\equiv
u+\varepsilon \varphi $ for $\varepsilon \in R.$ We computed the first
variation of $u$ in the direction $\varphi $ and obtained%
\begin{equation}
\frac{d\mathcal{F}(u_{\varepsilon })}{d\varepsilon }\mid _{\varepsilon =0\pm
}=\pm \int_{S(u)}\mid \nabla \varphi \mid +\int_{\Omega \backslash
S(u)}N^{u}\cdot \nabla \varphi +\int_{\Omega }H\varphi  \label{1.1}
\end{equation}

\noindent (see (3.3) in \cite{CHY}) where $S(u)$ denotes the singular set of 
$u,$ consisting of the points (called singular points) where $\nabla u+\vec{F%
}$ $=$ $0,$ and $N^{u}$ $\equiv $ $N_{\vec{F}}^{u}$ $\equiv $ $\frac{\nabla
u+\vec{F}}{|\nabla u+\vec{F}|}$ (called Legendrian or horizontal normal).
For a general $\vec{F}$ we cannot ignore the contribution of the first
integral in the right-hand side of (\ref{1.1}), caused by the singular set $%
S(u).$ For instance, we have $S(u)$ $=$ $\Omega $ in the case that $\vec{F}$ 
$=$ $0$ and $u$ $=$ $0.$ For the study of the size of singular set, we refer
the reader to \cite{Ba} (in which singular set is called characteristic
set). Let us denote the space of weakly differentiable functions by $%
W^{1}(\Omega )$ (see \cite{GT} for the precise definition)$.$ By (\ref{1.1})
we define weak solutions as follows:

\bigskip

\textbf{Definition 1.1.} (Definition 3.2 in \cite{CHY}) Let $\Omega \subset
R^{m}$ be a bounded domain. Let $\vec{F}$ be an $L_{loc}^{1}$ vector field
and $H$ be an $L_{loc}^{1}$ function on $\Omega .$ We say $u\in W^{1}(\Omega
)$ is a weak solution to equation $divN^{u}$ $=$ $H$ in $\Omega $ if for any 
$\varphi \in C_{0}^{\infty }(\Omega ),$ there holds%
\begin{equation}
\int_{S(u)}|\nabla \varphi |+\int_{\Omega \backslash S(u)}N^{u}\cdot \nabla
\varphi +\int_{\Omega }H\varphi \geq 0.  \label{1.1'}
\end{equation}

\noindent In this situation, we also say that $divN^{u}$ $=$ $H$ in the weak
sense.

\bigskip

With $\varphi $ replaced by $-\varphi $ in (\ref{1.1'}), it follows that
when the $m$-dimensional Hausdorff measure of $S(u),$ denoted as $\mathcal{H}%
^{m}(S(u)),$ vanishes, then $u$ $\in $ $W^{1}(\Omega )$ is a weak solution
to equation $divN^{u}$ $=$ $H$ in $\Omega $ if for any $\varphi \in
C_{0}^{\infty }(\Omega ),$ there holds

\begin{equation}
\int_{\Omega }N^{u}\cdot \nabla \varphi +\int_{\Omega }H\varphi =0.
\label{1.2}
\end{equation}%
\noindent That is, if $\mathcal{H}^{m}(S(u))$ $=$ $0,$ then $divN^{u}$ $=$ $%
H $ in the weak sense if (\ref{1.2}) holds. In this paper, we consider only
the situation that $S(u)$ is an empty set. A domain $\Omega $ is called
nonsingular (for $u$) if $S(u)$ is empty. So we can use (\ref{1.2}) as the
definition of weak solutions for a nonsingular domain $\Omega $. In this
paper we assume $u$ to be $C^1$, an assumption that is almost justified
according to a structure theorem of Franchi, Serapioni, and Serra Cassano (%
\cite{FSS01}), modulo a measure zero set, the major part of a very general
surface is $\boldsymbol{H}$-regular.

\bigskip

\textbf{Theorem 1.1} (\cite{FSS01}) \textit{If }$E$\textit{\ }$\subseteq $%
\textit{\ }$\boldsymbol{H}^{\boldsymbol{1}}$\textit{\ is a }$\boldsymbol{H}$%
\textit{-Caccioppoli set, then }$\partial _{\boldsymbol{H}}^{\ast }E,$%
\textit{\ the reduced boundary of }$E,$\textit{\ is }$\boldsymbol{H}$\textit{%
-rectifiable, that is}%
\begin{equation*}
\partial _{\boldsymbol{H}}^{\ast }E=\mathcal{N}\cup
\dbigcup\limits_{h=1}^{\infty }K_{h}
\end{equation*}%
\textit{\ \noindent where }$\mathcal{H}_{d}^{3}(\mathcal{N)}$\textit{\ }$=$%
\textit{\ }$0$\textit{\ and }$K_{h}$\textit{\ is a compact subset of a }$%
\boldsymbol{H}$\textit{-regular surface.}

\bigskip

On the other hand, it does not follow that a $\boldsymbol{H}$-regular
surface can be represented by a $C^1$ graph. We thank F. Serra-Cassano and
the referee for pointing out this fact.

Let $u$ $\in $ $C^{1}(\Omega ),$ $\vec{F}$ $\in $ $C^{0}(\Omega ),$ and $%
\Omega $ be nonsingular$.$ Then $N^{u}$ and $N^{u,\perp }$ exist and are in $%
C^{0}(\Omega ),$ where $N^{u,\perp }$ denotes the rotation of $N^{u}$ by $-%
\frac{\pi }{2}$ degrees$.$ By the O.D.E. theory (\cite{Ha}), the integral
curves of $N^{u,\perp }$ ($N^{u},$ resp.) exist and we call them
characteristic (seed, resp.) curves. In \cite{Pauls}, Pauls proved that for
a $C^{1}$ smooth (weak) solution to equation $divN^{u}$ $=$ $H$ with $\vec{F}
$ $=$ $(-y,$ $x)$ and $H$ $=$ $0$ in a nonsingular domain$,$ the
characteristic curves are straight lines and the seed curves are $C^{2}$
smooth under the condition that $N^{u}$ $\in $ $W^{1,1}.$ We will show that
the condition: $N^{u}$ $\in $ $W^{1,1}$ is not necessary while the curvature
(along $N^{u,\perp }$ direction) of a characteristic curve is $-H$ under
mild regularity condition on $H$ (and nothing to do with precise form of $%
\vec{F};$ see Theorem A, Theorem D below)

For $u$ $\in $ $C^{1}(\Omega ),$ $\vec{F}$ $\in $ $C^{0}(\Omega )$ in a
nonsingular plane domain $\Omega ,$ since $|N^{u}|$ $=$ $1,$ we can write $%
N^{u}$ $=$ $(\cos \theta ,$ $\sin \theta )$ locally with $\theta $ $\in $ $%
C^{0}.$ We may forget $u$ and consider $\theta $ $\in $ $C^{0}$ locally as
an independent variable, then define $N$ $\equiv $ $(\cos \theta ,$ $\sin
\theta )$ such that $N$ and $N^{\perp }$ $\equiv $ $(\sin \theta ,$ $-\cos
\theta )$ are $C^{0}$ vector fields. We also call the integral curves of $%
N^{\perp }$ ($N,$ resp.) characteristic (seed, resp.) curves (see (\ref{1.5}%
) below). Similarly to (\ref{1.2}) for $\theta $ $\in $ $C^{0}(\Omega )$ we
define 
\begin{equation}
\func{div}N\equiv \func{div}(\cos \theta ,\sin \theta )\equiv (\cos \theta
)_{x}+(\sin \theta )_{y}=H  \label{1.3}
\end{equation}

\noindent in the weak sense, meaning 
\begin{equation*}
\int_{\Omega }N\cdot \nabla \varphi +\int_{\Omega }H\varphi =0.
\end{equation*}%
\noindent for any $\varphi \in C_{0}^{\infty }(\Omega )$. On the other hand,
equation $divN^{u}$ $=$ $H$ (i.e., $\theta $ arises from $u)$ provides more
information. For $u$ $\in $ $C^{1},$ $\vec{F}$ $\in $ $C^{0},$ let $D$ $%
\equiv $ $|\nabla u+\vec{F}|.$ If $D$ $\neq $ $0$ in $\Omega ,$ write $%
N^{u,\perp }$ $=$ $(\sin \theta ,$ $-\cos \theta )$ and $\vec{F}^{\perp }$ $%
= $ $(F_{2},$ $-F_{1})$ for $\vec{F}$ $=$ $(F_{1},$ $F_{2}).$ Then for any
Lipschitzian domain $\Omega ^{\prime }$ $\subset \subset $ $\Omega ,$ we have%
\begin{equation}
\int_{\Omega ^{\prime }}(DN^{u,\perp }-\vec{F}^{\perp })\cdot \nabla \varphi
=\int_{\Omega ^{\prime }}(u_{y},-u_{x})\cdot \nabla \varphi =0  \label{1.4}
\end{equation}

\noindent (see the proof of Lemma 3.1) The \textquotedblright
integrability\textquotedblright\ condition (\ref{1.4}) (due to
\textquotedblright $u_{yx}$ $=$ $u_{xy}"$) makes equation $divN^{u}$ $=$ $H$
have more properties than (\ref{1.3}).

Now we study (\ref{1.3}) as a single equation for $\theta .$ For $\theta $ $%
\in $ $C^{1}$, (\ref{1.3}) is a first-order equation whose characteristic
curve $\Gamma $ $=$ \{$(x(\sigma ),$ $y(\sigma ))\}$ $\in $ $C^{1}$
satisfies 
\begin{eqnarray}
\frac{dx}{d\sigma } &=&\sin \theta (x(\sigma ),y(\sigma )),  \label{1.5} \\
\frac{dy}{d\sigma } &=&-\cos \theta (x(\sigma ),y(\sigma )).  \notag
\end{eqnarray}

\noindent Note that $\sigma $ is a unit-speed parameter of $\Gamma .$ For $%
\theta $ $\in $ $C^{0},$ we still use (\ref{1.5}) as the definition of
characteristic curves.

\bigskip

\textbf{Theorem A. }\textit{Let }$\Omega $\textit{\ be a domain of }$R^{2}$%
\textit{\ and }$H\in C^{0}(\Omega ).$\textit{\ Let }$\theta $ $\in $ $%
C^{0}(\Omega )$\textit{\ satisfy equation (\ref{1.3}) in the weak sense,
i.e. }%
\begin{equation}
\int_{\Omega }(\cos \theta ,\sin \theta )\cdot \nabla \varphi +\int_{\Omega
}H\varphi =0  \label{1.6}
\end{equation}

\noindent \textit{for all }$\varphi \in C_{0}^{\infty }(\Omega )$\textit{\
(cf. (\ref{1.2}))}$.$\textit{\ Let }$\Gamma \subset \Omega $\textit{\ be a (}%
$C^{1}$ \textit{smooth)} \textit{characteristic curve with }$\sigma $ 
\textit{being the unit-speed parameter, satisfying (\ref{1.5}). Then }$%
\Gamma $\textit{\ is }$C^{2}$\textit{\ smooth and the curvature of }$\Gamma $%
\textit{\ (along }$N^{\perp }$ \textit{direction}) \textit{equals }$-H,$ 
\textit{that is, }$\frac{d\theta }{d\sigma }=-H.$

\bigskip

Note that $\theta $ $\in $ $C^{0}$ implies $N$ $\equiv $ $(\cos \theta ,$ $%
\sin \theta )$ $\in $ $C^{0}$ and $N$ $=$ $N^{u}$ $\equiv $ $N_{\vec{F}%
}^{u}\in $ $C^{0}$ if it arises from $u$ $\in $ $C^{1}$ and $\vec{F}$ $\in $ 
$C^{0}.$ Recall that in Theorem A of \cite{Pauls}, Pauls considered the $H$ $%
=$ $0$ case, in which $\Gamma $ is a straight line under the condition that
components of the horizontal Gauss map (i.e., $N_{\vec{F}}^{u}$ in our
notation with $\vec{F}$ $=$ $(-y,$ $x)$) are in $W^{1,1}(\Omega ).$ In
Theorem A above, if $\theta $ satisfies (\ref{1.6}), we prove that $\Gamma $
is a minimizer for the following energy functional:%
\begin{equation*}
L_{H}(\Gamma )\equiv |\Gamma |-\int_{\Omega _{\Gamma }}Hdxdy
\end{equation*}%
\noindent where $|\Gamma |$ denotes the length of $\Gamma $ (see Section 2
for the definition of $\Omega _{\Gamma }).$ So the basic Calculus of
Variations tells us that the curvature of $\Gamma $\ (along $N^{\perp }$
direction) equals $-H$ without invoking extra regularity assumption$.$ Also $%
H$ is only required to be $C^{0}.$ In \cite{MR} Monti and Rickly considered
the case of $H$ $=$ constant for a convex isoperimetric set. We do not need
convexity in Theorem A.

For $\theta $ $\in $ $C^{0},$ $N$ ($N^{\perp },$ resp.) $\equiv $ $(\cos
\theta ,$ $\sin \theta )$ ($\equiv $ $(\sin \theta ,$ $-\cos \theta ),$
resp.) is a $C^{0}$ vector field. Then for any $p$ $\in $ $\Omega ,$ there
exists at least one integral curve, i.e., seed curve (characteristic curve,
resp.) passing through $p$. The uniqueness of integral curves for a $C^{0}$
vector field does not hold in general (see page 18 in \cite{Ha}). In Section
3 we will prove uniqueness theorems B and B$^{\prime }$ for characteristic
and seed curves (see below). Let $p$ $\in $ $\Omega $ and $B_{r}(p)$ $\equiv 
$ $\{q$ $\in $ $\Omega $ $|$ $|q-p|$ $<$ $r\}.$ Define $H_{M}(r)$ $\equiv $ $%
\max_{q\in \partial B_{r}(p)}$ $|H(q)|.$

\bigskip

\textbf{Theorem B.} \textit{(a)\ Let }$\theta $\textit{\ }$\in $\textit{\ }$%
C^{0}(\Omega )$\textit{\ and }$H$\textit{\ }$\in $\textit{\ }$%
L_{loc}^{1}(\Omega )$\textit{\ satisfy (\ref{1.6})}$.$\textit{\ Let }$p$%
\textit{\ }$\in $\textit{\ }$\Omega $\textit{\ and suppose there is }$r_{0}$%
\textit{\ }$>$\textit{\ }$0$\textit{\ such that }$B_{r_{0}}(p)$\textit{\ }$%
\subset \subset $\textit{\ }$\Omega $\textit{\ and }%
\begin{equation}
\int_{0}^{r_{0}}H_{M}(r)dr<\infty .  \label{1.7}
\end{equation}%
\textit{Then\ there is }$r_{1},$\textit{\ }$0$\textit{\ }$<$\textit{\ }$%
r_{1} $\textit{\ }$\leq $\textit{\ }$r_{0},$\textit{\ such that there exists
a unique seed curve passing through }$p$\textit{\ in }$B_{r_{1}}(p).$\textit{%
\ }

\textit{(b) Let }$\theta $\textit{\ }$\in $\textit{\ }$C^{0}(\Omega )$%
\textit{\ and }$H$\textit{\ }$\in $\textit{\ }$C^{0,1}(\Omega )$\textit{\
(Lipschitzian) satisfy (\ref{1.6}). Then for any point }$p$\textit{\ }$\in $%
\textit{\ }$\Omega ,$\textit{\ we can find }$r_{1}$ $>$ $0$ \textit{such
that there exists a unique characteristic curve passing through }$p$\textit{%
\ in }$B_{r_{1}}(p).$

\textit{\bigskip }

In Theorem B (b), if $H$ is only continuous, we can give an example for the
nonuniqueness of characteristic curves (see Example 3.2). Note that $u$ is
not involved in Theorems A and B. Now we consider $u.$ Let $u$ $\in $ $C^{1}$
and $\vec{F}$ $=$ $(F_{1},$ $F_{2})$ $\in $ $C^{1}.$ Recall that a point $p$ 
$\in $ $\Omega $ $\subset $ $R^{2}$ is called singular (nonsingular, resp.)
if $\nabla u$ $+$ $\vec{F}$ $=$ $0$ ($\neq $ $0,$ resp.$)$ at $p.$ At a
nonsingular point, we recall $N$ $\equiv $ $N^{u}$ $\equiv $ $\frac{\nabla u+%
\vec{F}}{|\nabla u+\vec{F}|}.$ We call $\Omega $ nonsingular if every point
of $\Omega $ is not singular. We have another uniqueness theorem for
characteristic curves.

\bigskip

\textbf{Theorem B}$^{\prime }$. \textit{Let }$u$\textit{\ }$:$\textit{\ }$%
\Omega $\textit{\ }$\subset $\textit{\ }$R^{2}$\textit{\ }$\rightarrow $%
\textit{\ }$R$\textit{\ be a }$C^{1}$\textit{\ smooth function such that }$%
\Omega $\textit{\ is a nonsingular domain with }$\vec{F}$ $\in $ $%
C^{1}(\Omega )$\textit{. Then for any point }$p$\textit{\ }$\in $\textit{\ }$%
\Omega ,$\textit{\ we can find }$r_{1}$ $>$ $0$ \textit{such that }$%
B_{r_{1}}(p)$\textit{\ }$\subset $\textit{\ }$\Omega $\textit{\ and there
exists a unique characteristic curve passing through }$p$\textit{\ in }$%
B_{r_{1}}(p).$

\bigskip

In Theorem B$^{\prime }$ we only assume $u,$ $\vec{F}$ $\in $ $C^{1},$ and
do not use any property of $H,$ in contrast to Theorem B (b)$.$ Even for the
case $H$ $=$ $0,$ seed curves may only be $C^{1}$ smooth, but not $C^{2}$
smooth (see the remark after the proof of Theorem D in Section 5). However
if $N$ $\equiv $ $(\cos \theta ,$ $\sin \theta )$ arises from $u$ (i.e., $N$ 
$=$ $N^{u};$ see (\ref{1.2})), Pauls (\cite{Pauls}) proved that when $u$ $%
\in $ $C^{1},$ $\theta $ $\in $ $C^{0}\cap W^{1,1},$ and $H$ $=$ $0,$ then
the seed curves are $C^{2}$ smooth. In Theorem D we prove the same
conclusion under the condition that $u$ $\in $ $C^{1}$ ($\theta $ $\in $ $%
C^{0}$ follows) and $H$ $\in $ $C^{1}$ (in fact, that $H$ $\in $ $C^{0}$ and
only $C^{1}$ in the $N$ direction is enough)$.$

If $u$ $\in $ $C^{1}(\Omega )$ also satisfies (\ref{1.2}) with $\vec{F}$ $%
\in $ $C^{1}(\Omega )$ and $S(u)$ is empty in $\Omega ,$ can we have higher
order regularity for $u,$ say, $u$ $\in $ $C^{2}?$ This is impossible as
shown by the following example. Let $u_{g}$ $\equiv $ $xy$ $+$ $g(y)$ where $%
g$ $\in $ $C^{1}\backslash C^{2}.$ Then $u_{g}$ satisfies (\ref{1.2}) with $%
H $ $=$ $0,$ $\vec{F}$ $=$ $(-y,x)$ on any nonsingular domain $\Omega $ for $%
u_{g}.$ On the other hand, the characteristic and seed curves associated to $%
u_{g}$ are all the same for different $g$'s. That is, $g$ determines the
differentiability of $u_{g},$ but does not affect the shape of
characteristic and seed curves. So we can prove that $\theta $ is in fact $%
C^{1}$ smooth (hence $N$ $\in $ $C^{1},$ but not $u$ $\in $ $C^{2})$ (see
Theorem D below)$.$ Before doing this we need to introduce some kind of
special coordinates.

\bigskip

\textbf{Definition 1.2.} Let $N$ be a $C^{0}$ vector field with $|N|$ $%
\equiv $ $1$ on a domain $\Omega $ $\subset $ $R^{2}.$ A system of $C^{1}$
smooth local coordinates $s,$ $t$ is called a system of characteristic
coordinates if $s$ and $t$ have the property that $\nabla s$ $\parallel $ $%
N^{\perp }$ and $\nabla t$ $\parallel $ $N,$ i.e., $\nabla s$ and $\nabla t$
are parallel to $N^{\perp }$ and $N$, resp.$.$ It follows that $\{t$ $=$
constant$\}$ are characteristic curves while $\{s$ $=$ constant$\}$ are seed
curves.

\bigskip

Let $\Gamma _{(x,y)}(\sigma )$ ($\Lambda _{(x,y)}(\tau ),$ resp.) denote a
characteristic (seed, resp.) curve passing through $(x,y),$ parametrized by
the arc length $\sigma $ ($\tau ,$ resp.) with $\frac{d\Gamma
_{(x,y)}(\sigma )}{d\sigma }$ $=$ $N^{\perp }$ ($\frac{d\Lambda
_{(x,y)}(\tau )}{d\tau }$ $=$ $N,$ resp.)$.$ For continuous functions $g,$ $%
f $ we write 
\begin{equation*}
(N^{\perp }g)(x,y)\equiv \frac{dg(\Gamma _{(x,y)}(\sigma ))}{d\sigma }\mid
_{\sigma =\sigma _{0}}
\end{equation*}

\noindent $((Nf)(x,y)\equiv \frac{df(\Lambda _{(x,y)}(\tau ))}{d\tau }\mid
_{\tau =\tau _{0}},$ resp.) if exists, where $(x,$ $y)$ $=$ $\Gamma
_{(x,y)}(\sigma _{0})$ ($(x,$ $y)$ $=$ $\Lambda _{(x,y)}(\tau _{0}),$ resp.)$%
.$ For a planar $C^{1}$ vector field $\vec{F}$ $=$ $(F_{1},$ $F_{2}),$ we
define%
\begin{equation*}
rot\vec{F}:=(F_{2})_{x}-(F_{1})_{y}.
\end{equation*}

\noindent We construct a system of characteristic coordinates in the
following theorem.

\bigskip

\textbf{Theorem C.} \textit{Let }$u$\textit{\ }$:$\textit{\ }$\Omega $%
\textit{\ }$\subset $\textit{\ }$R^{2}$\textit{\ }$\rightarrow $\textit{\ }$%
R $\textit{\ be a }$C^{1}$\textit{\ smooth solution to (\ref{1.2}) (}$\Omega 
$\textit{\ being nonsingular) with }$\vec{F}$ $\in $ $C^{1}(\Omega )$ 
\textit{and} $H$\textit{\ }$\in $\textit{\ }$C^{0}(\Omega ).$\textit{\ Then
for any point }$p_{0}$\textit{\ }$\in $\textit{\ }$\Omega $\textit{\ there
exist a neighborhood }$\Omega ^{\prime }$\textit{\ }$\subset $\textit{\ }$%
\Omega $\textit{\ and real functions }$s,$\textit{\ }$t$\textit{\ }$\in $%
\textit{\ }$C^{1}(\Omega ^{\prime })$\textit{\ such that }$\{t$\textit{\ }$=$%
\textit{\ constants}$\}$\textit{\ and }$\{s$\textit{\ }$=$\textit{\ constants%
}$\}$\textit{\ are characteristic curves and seed curves, respectively.
Moreover, there are positive functions }$f,$\textit{\ }$g$\textit{\ }$\in $%
\textit{\ }$C^{0}(\Omega ^{\prime })$\textit{\ such that }%
\begin{equation}
\nabla s=fN^{\perp },\nabla t=gDN.  \label{1.8}
\end{equation}

\noindent \textit{Also }$Nf$\textit{\ and }$N^{\perp }g$\textit{\ exist and
are continuous in }$\Omega ^{\prime }$\textit{. In fact, }$f$\textit{\ and }$%
g$\textit{\ satisfy the following equations}%
\begin{equation}
Nf+fH=0,N^{\perp }g+\frac{(rot\vec{F})g}{D}=0.  \label{1.9}
\end{equation}

\noindent \textit{For a perhaps smaller neighborhood }$\Omega ^{\prime
\prime }$\textit{\ }$\subset $\textit{\ }$\Omega ^{\prime }$\textit{\ of }$%
p_{0},$\textit{\ the map }$\Psi $\textit{\ }$:$\textit{\ }$(x,$\textit{\ }$%
y) $\textit{\ }$\in $\textit{\ }$\Omega ^{\prime \prime }$\textit{\ }$%
\rightarrow $\textit{\ }$(s,$\textit{\ }$t)$\textit{\ }$\in $\textit{\ }$%
\Psi (\Omega ^{\prime \prime })$\textit{\ is a }$C^{1}$\textit{\
diffeomorphism such that}%
\begin{equation}
\Psi ^{\ast }(\frac{ds^{2}}{f^{2}}+\frac{dt^{2}}{g^{2}D^{2}})=dx^{2}+dy^{2}.
\label{1.10}
\end{equation}

\bigskip

We remark that the existence of $C^{1}$ smooth $s$ can be proved for $N$
satisfying (\ref{1.6}) (i.e., not defined by $u)$ instead of (\ref{1.2})
(see Theorem 4.1).

Recall that for $u$ $\in $ $C^{1}(\Omega ),$ $\vec{F}$ $\in $ $C^{0}(\Omega
),$ and $\Omega $ being nonsingular, there exists $\theta $ $\in $ $C^{0}$
locally such that $N^{u}$ $=$ $(\cos \theta ,$ $\sin \theta ).$

\bigskip

\textbf{Corollary C.1}. \textit{Suppose we are in the situation of Theorem
C. Then }$\theta $\textit{\ is }$C^{1}$\textit{\ smooth in }$s$\textit{\ and
there holds}%
\begin{equation}
\theta _{s}\equiv \frac{\partial \theta }{\partial s}=-\frac{H}{f}.
\label{1.10.1}
\end{equation}

\bigskip

Note that by $\theta $ being $C^{1}$ smooth in $s,$ we mean that $\theta
_{s} $ $\equiv $ $\frac{\partial \theta }{\partial s}$ exists and is
continuous. Since $f$ is $C^{1}$ smooth in $t,$ $\theta _{s}$ is also $C^{1}$
smooth in $t$ if we assume that $H$ has the same property according to (\ref%
{1.10.1}). In fact, we can prove that $\theta $ is $C^{1}$ smooth in $t$
too, and hence $\theta $ $\in $ $C^{1}.$ That is, $\theta $ gains one
derivative.

\bigskip

\textbf{Theorem D.} \textit{Let }$u$\textit{\ }$:$\textit{\ }$\Omega $%
\textit{\ }$\subset $\textit{\ }$R^{2}$\textit{\ }$\rightarrow $\textit{\ }$%
R $\textit{\ be a }$C^{1}$\textit{\ smooth (weak) solution to }$divN^{u}$ $=$
$H$\textit{\ in }$\Omega $ \textit{(}$\Omega $\textit{\ being nonsingular)
with }$\vec{F}$ $\in $ $C^{1}(\Omega )$ \textit{and }$H$\textit{\ }$\in $%
\textit{\ }$C^{0}(\Omega ).$\textit{\ Suppose }$N(H)$\textit{\ exists and is
continuous.\ Then }$\theta $\textit{\ }$\in $\textit{\ }$C^{1}$ \textit{and
the characteristic and seed curves are }$C^{2}$\textit{\ smooth}$.$\textit{\
Moreover, }$N^{\perp }D$\textit{\ exists and is continuous in }$\Omega $%
\textit{. In }$(s,t)$\textit{\ coordinates near a given point as in Theorem
C, we have}%
\begin{equation}
\theta _{t}\equiv \frac{\partial \theta }{\partial t}=\frac{rot\vec{F}}{%
gD^{2}}-\frac{N^{\perp }(\log D)}{gD}=\frac{1}{gD^{2}}(rot\vec{F}-N^{\perp
}D).  \label{1.11}
\end{equation}

\textit{\bigskip }

We remark that in case $H$ $=$ $0$ or constant, we can prove Theorem D
directly from the precise parametric expression of $x$ or $y.$ The situation 
$H$ $=$ a nonzero constant arises from considering the boundary of a $C^{2}$
isoperimetric set. Pansu (\cite{Pan82}) conjectured that an isoperimetric
set is congruent with a certain type of sphere. In \cite{RR}, Ritor\'{e} and
Rosales proved Pansu's conjecture for isoperimetric sets of class $C^{2}$
without any symmetry assumption. Later Monti and Rickly (\cite{MR}) obtained
the same result for convex isoperimetric sets without regularity assumptions.

We outline the proof of $\theta $\textit{\ }$\in $\textit{\ }$C^{1}$ in
Theorem D as follows. According to (\ref{1.10.1}), we have good control for $%
\theta $ along the characteristic curves, i.e. the $s$-direction. If the
control for $\theta $ fails along a seed curve, i.e. $t$-direction, say, at
some $s_{0},$ then we show that it fails also for $s$ near $s_{0}.$ That is,
the jump of a certain concerned quantity is kept in short \textquotedblright 
$s$-time\textquotedblright\ along the characteristic curves. This ends up to
reach a contradiction. We borrow the idea of conveying information along the
characteristic curves from the study of hyperbolic P.D.E. (\cite{John}).
After this paper was submitted, we were informed of recent results about
regularity (\cite{CCM1}, \cite{CCM2}) by the referee. In \cite{CCM2},
Capogna, Citti, and Manfredini proved an interesting result, among others,
that the Lipschitz minimizer obtained in Theorem A of \cite{CHY} is actually 
$C^{1,\alpha }$ in a neighborhood of a nonsingular point under some extra
condition. We would like to thank the referee for useful information,
detailed comments, and grammatical suggestions.

\bigskip

\section{\textbf{Curvature of characteristic curves-proof of Theorem A}}

The following lemma should be a standard result. For completeness we give a
proof.

\bigskip

\textbf{Lemma 2.1}. \textit{Let }$\theta $\textit{\ }$\in $\textit{\ }$%
C^{0}(\Omega )$\textit{\ and }$H$\textit{\ }$\in $\textit{\ }$%
L_{loc}^{1}(\Omega )$\textit{\ satisfy (\ref{1.6}). Then for any
Lipschitzian domain }$\Omega ^{\prime }$\textit{\ }$\subset \subset $\textit{%
\ }$\Omega $\textit{\ and }$\varphi $\textit{\ }$\in $\textit{\ }$%
C^{1}(\Omega ),$\textit{\ there holds}%
\begin{equation}
\oint_{\partial \Omega ^{\prime }}\varphi N\cdot \nu =\int_{\Omega ^{\prime
}}(\nabla \varphi )\cdot N+\varphi H  \label{2.1}
\end{equation}%
\textit{\noindent where }$N\equiv (\cos \theta ,\sin \theta )$ \textit{and }$%
\nu $\textit{\ denotes the unit outer normal to }$\partial \Omega ^{\prime
}. $

\bigskip

\proof
Take a Lipschitzian domain $\Omega ^{\prime \prime }$ such that $\Omega
^{\prime }$ $\subset \subset $ $\Omega ^{\prime \prime }$ $\subset \subset $ 
$\Omega .$ Then there exists a sufficiently small number $\varepsilon _{0}$ $%
>$ $0$ such that for any $\psi $ $\in $ $C_{0}^{1}(\Omega ^{\prime \prime })$
and any $\varepsilon ,$ $0$ $<$ $\varepsilon $ $<$ $\varepsilon _{0},$ the
mollifier $\psi _{\varepsilon }$ $\in $ $C_{0}^{\infty }(\Omega ).$ From (%
\ref{1.6}) we have

\begin{equation*}
0=\int_{\Omega }N\cdot \nabla \psi _{\varepsilon }+H\psi _{\varepsilon
}=\int_{\Omega ^{\prime \prime }}N_{\varepsilon }\cdot \nabla \psi
+H_{\varepsilon }\psi .
\end{equation*}

\noindent It then follows that $div$ $N_{\varepsilon }$ $=$ $H_{\varepsilon
} $ (strong sense) in $\Omega ^{\prime \prime }$ and since $\varphi $ $\in $ 
$C^{1}(\Omega ),$ we have%
\begin{equation}
\oint_{\partial \Omega ^{\prime }}\varphi N_{\varepsilon }\cdot \nu
=\int_{\Omega ^{\prime }}(\nabla \varphi )\cdot N_{\varepsilon }+\varphi
H_{\varepsilon }.  \label{2.2}
\end{equation}

\noindent Since $N$ $\in $ $C^{0}(\Omega ),$ $N_{\varepsilon }$ converges to 
$N$ uniformly in $\Omega ^{\prime \prime }$ while $H_{\varepsilon }$
converges to $H$ in $L^{1}(\Omega ^{\prime \prime })$ as $\varepsilon
\rightarrow 0.$ Therefore letting $\varepsilon \rightarrow 0$ in (\ref{2.2})
gives (\ref{2.1}).

\bigskip 
\endproof%

By taking $\varphi $ $\equiv $ $1$ in (\ref{2.1}), we obtain

\bigskip

\textbf{Corollary 2.2}. \textit{Let }$\theta $\textit{\ }$\in $\textit{\ }$%
C^{0}(\Omega )$\textit{\ and }$H$\textit{\ }$\in $\textit{\ }$%
L_{loc}^{1}(\Omega )$\textit{\ satisfy (\ref{1.6}). Then for any
Lipschitzian domain }$\Omega ^{\prime }$\textit{\ }$\subset \subset $\textit{%
\ }$\Omega ,$\textit{\ we have}%
\begin{equation}
\oint_{\partial \Omega ^{\prime }}N\cdot \nu =\int_{\Omega ^{\prime }}H
\label{2.3}
\end{equation}

\noindent \textit{where }$N\equiv (\cos \theta ,\sin \theta ).$

\bigskip

\proof%
{\large \ (of Theorem A)}{\normalsize . }Without loss of generality, we may
assume that locally $\Gamma $ is a piece of a $C^{1}$ smooth graph $(x,y(x))$
with $0$ $\leq $ $x$ $\leq $ $a$ ($a$ $>$ $0$) where $y(x)$ $>$ $0$, $%
y^{\prime }(x)$ is bounded, and the domain $\Omega _{\Gamma }$ surrounded by 
$\Gamma $ and the three line segments connecting $(0,y(0))$, $(0,0),$ $%
(a,0), $ $(a,y(a))$ satisfies $\Omega _{\Gamma }$ $\subset \subset $ $\Omega
.$ We may also assume that $N^{\perp }$ $\equiv $ $(\sin \theta ,$ $-\cos
\theta )$ $=$ $\frac{(1,y^{\prime })}{\sqrt{1+(y^{\prime })^{2}}}$ on $%
\Gamma .$ Let $\Gamma _{\varepsilon }$ be a family of small perturbations of 
$\Gamma ,$ having the same endpoints $(0,y(0)),$ $(a,y(a))$ and described by 
$(x,y(x)+\varepsilon \varphi (x))$ ($\varphi $ $\in $ $C_{0}^{\infty
}([0,a])) $ (see Figure 2.1)

\begin{figure}[htp]
\centering
\begin{pspicture}(-1,-1.5)(8,5) 
\psset{algebraic=true, xunit=1.2,yunit=0.8, linewidth=0.5pt} 
\psaxes[ticks=none, labels=none]{->}(0,0)(-1,-1)(7,5) 
\psecurve[linewidth=1.0pt]{-}(-0.5,1.8)(0,2.5)(1.2,3)(2.6,2.7)(4,3)(5,3.4)(6,3.5)(6.8,3.4)
\psecurve[linewidth=1.0pt]{-}(-0.5,1.8)(0,2.5)(1.2,3.5)(2.6,3.2)(4,3.5)(5,3.9)(6,3.5)(6.8,3.4)
\psline[linewidth=1.0pt]{-}(6,0)(6,3.5)
\uput{1pt}[90]{0}(3,3.6){$\Gamma_\epsilon$}
\uput{1pt}[-90]{0}(3,2.4){$\Gamma$}
\uput{3pt}[90]{0}(0,5){$y$}
\uput{3pt}[0]{0}(7,0){$x$}
\uput{6pt}[225]{0}(0,0){$(0,0)$}
\uput{6pt}[-45]{0}(6,0){$(a,0)$}
\uput{3pt}[180]{0}(0,2.5){$(0,y(0))$}
\uput{3pt}[0]{0}(6,3.5){$(a,y(a))$}

\uput*{4pt}[-90]{0}(3.2,-1.5){Figure 2.1} 
\end{pspicture}
\end{figure}

The domain $\Omega _{\Gamma _{\varepsilon }}$ are defined similarly and $%
\Omega _{\Gamma _{\varepsilon }}$ $\subset \subset $ $\Omega $ for $%
|\varepsilon |$ small enough. Note that $y_{\varepsilon }(x)$ $\equiv $ $%
y(x) $ $+$ $\varepsilon \varphi (x)$ $>$ $0$ (on $[0,a]$) also for $%
|\varepsilon | $ small enough$.$ Let $\sigma $ denote the arc length
parameter. Let%
\begin{equation*}
G(\Gamma _{\varepsilon })\equiv \int_{\Gamma _{\varepsilon }}N\cdot \nu
d\sigma -\int_{\Omega _{\Gamma _{\varepsilon }}}Hdxdy
\end{equation*}

\noindent where $\nu $ denotes the unit outer normal of $\Omega _{\Gamma
_{\varepsilon }}$. Observe that

\begin{eqnarray}
&&G(\Gamma _{\varepsilon })-G(\Gamma )  \label{2.4} \\
&=&\{\oint_{\partial \Omega _{\Gamma _{\varepsilon }}}N\cdot \nu d\sigma
-\int_{\Omega _{\Gamma _{\varepsilon }}}Hdxdy\}-\{\oint_{\partial \Omega
_{\Gamma }}N\cdot \nu d\sigma -\int_{\Omega _{\Gamma }}Hdxdy\}  \notag \\
&=&0-0=0  \notag
\end{eqnarray}

\noindent by (\ref{2.3}). Let $|\Gamma |$ denote the length of $\Gamma .$
Along $\Gamma $ we have $\nu $ $=$ $N$ since $N^{\perp }$ is tangent to $%
\Gamma ,$ a characteristic curve, $y(x)$ $>$ $0$ and $N^{\perp }$ $=$ $\frac{%
(1,y^{\prime })}{\sqrt{1+(y^{\prime })^{2}}}$ on $\Gamma $. It follows that

\begin{eqnarray}
|\Gamma |-\int_{\Omega _{\Gamma }}Hdxdy &=&G(\Gamma )=G(\Gamma _{\varepsilon
})\text{ (by (\ref{2.4}))}  \label{2.5} \\
&\leq &\int_{\Gamma _{\varepsilon }}|N\cdot \nu |d\sigma -\int_{\Omega
_{\Gamma _{\varepsilon }}}Hdxdy  \notag \\
&\leq &|\Gamma _{\varepsilon }|-\int_{\Omega _{\Gamma _{\varepsilon }}}Hdxdy%
\text{ \ (by }|N\cdot \nu |\leq |N||\nu |=1).  \notag
\end{eqnarray}

\noindent Let 
\begin{equation*}
L_{H}(\Gamma )\equiv |\Gamma |-\int_{\Omega _{\Gamma }}Hdxdy.
\end{equation*}%
\noindent We learn from (\ref{2.5}) that $L_{H}(\Gamma )$ is the minimum of $%
L_{H}(\Gamma _{\varepsilon }).$ Therefore for $\varepsilon $ $\in $ $%
R\backslash \{0\},$ $|\varepsilon |$ small enough, we have

\begin{eqnarray}
0 &\leq &|\varepsilon |^{-1}\{(|\Gamma _{\varepsilon }|-\int_{\Omega
_{\Gamma _{\varepsilon }}}Hdxdy)-(|\Gamma |-\int_{\Omega _{\Gamma }}Hdxdy)\}
\label{2.6} \\
&=&|\varepsilon |^{-1}(|\Gamma _{\varepsilon }|-|\Gamma |)-|\varepsilon
|^{-1}(\int_{\Omega _{\Gamma _{\varepsilon }}}Hdxdy-\int_{\Omega _{\Gamma
}}Hdxdy).  \notag
\end{eqnarray}%
\noindent Compute

\begin{eqnarray}
&&\varepsilon ^{-1}(|\Gamma _{\varepsilon }|-|\Gamma |)  \label{2.6.1} \\
&=&\varepsilon ^{-1}\{\int_{0}^{a}\sqrt{1+(y^{\prime }(x)+\varepsilon
\varphi ^{\prime }(x))^{2}}dx-\int_{0}^{a}\sqrt{1+(y^{\prime }(x))^{2}}dx\} 
\notag \\
&=&\int_{0}^{a}\frac{2y^{\prime }(x)\varphi ^{\prime }(x)+\varepsilon
(\varphi ^{\prime }(x))^{2}}{\sqrt{1+(y^{\prime }(x)+\varepsilon \varphi
^{\prime }(x))^{2}}+\sqrt{1+(y^{\prime }(x))^{2}}}dx  \notag \\
&\rightarrow &\int_{0}^{a}\frac{y^{\prime }(x)}{\sqrt{1+(y^{\prime }(x))^{2}}%
}\varphi ^{\prime }(x)dx  \notag
\end{eqnarray}

\noindent while for some $\tilde{y}_{\varepsilon }(x)$ between $y(x)$ and $%
y(x)+\varepsilon \varphi (x)$ by the mean value theorem$,$ we have%
\begin{eqnarray}
&&\varepsilon ^{-1}(\int_{\Omega _{\Gamma _{\varepsilon
}}}Hdxdy-\int_{\Omega _{\Gamma }}Hdxdy)  \label{2.6.2} \\
&=&\varepsilon ^{-1}\int_{0}^{a}H(x,\tilde{y}_{\varepsilon }(x))\varepsilon
\varphi (x)dx  \notag \\
&\rightarrow &\int_{0}^{a}H(x,y(x))\varphi (x)dx  \notag
\end{eqnarray}

\noindent as $\varepsilon $ $\rightarrow $ $0$ by Lebesgue's dominated
convergence theorem.

From (\ref{2.6}), (\ref{2.6.1}), and (\ref{2.6.2}) we obtain 
\begin{equation}
\frac{d}{dx}(\frac{y^{\prime }(x)}{\sqrt{1+(y^{\prime }(x))^{2}}})=-H(x,y(x))
\label{2.7}
\end{equation}

\noindent in the weak sense. Since $\frac{y^{\prime }(x)}{\sqrt{1+(y^{\prime
}(x))^{2}}}$ and $H(x,y(x))$ are continuous in $x,$ we actually have $\frac{%
y^{\prime }(x)}{\sqrt{1+(y^{\prime }(x))^{2}}}$ $\in $ $C^{1}$ with respect
to $x$ and (\ref{2.7}) holds in the strong sense. Also it follows that $y$ $%
\in $ $C^{2}$ since $y^{\prime }(x)$ $=$ $\frac{h(x)}{\sqrt{1-(h(x))^{2}}}$ $%
\in $ $C^{1}$ where $h(x)$ $=$ $\frac{y^{\prime }(x)}{\sqrt{1+(y^{\prime
}(x))^{2}}}.$ We have proved that $\Gamma $\textit{\ }is $C^{2}$\ smooth.
Since $(\cos \theta ,$ $\sin \theta )$ $=$ $\frac{(1,y^{\prime })}{\sqrt{%
1+(y^{\prime })^{2}}},$ we have 
\begin{equation*}
\frac{d\sin \theta }{dx}=-H
\end{equation*}

\noindent by (\ref{2.7}). It follows from $\frac{dx}{d\sigma }$ $=$ $\cos
\theta $ that $-H$ $=$ $\cos \theta \frac{d\theta }{dx}$ $=$ $\frac{d\theta 
}{d\sigma }.$

\bigskip \textbf{\bigskip }%
\endproof%

\section{Uniqueness of Characteristic and seed curves}

Since $N$ and $N^{\perp }$ are $C^{0}$ vector fields, there exist integral
curves (called seed and characteristic curves, resp.) of $N$ and $N^{\perp }$
passing through any given point. Uniqueness does not hold in general (see 
\cite{Ha}, page 18). However, if $N$ satisfies (\ref{1.6}), we have
uniqueness.

\bigskip

\proof
({\large of Theorem B)} First we will prove (a). Since $N$ $\in $ $%
C^{0}(\Omega ),$ we can choose $r_{2},$ $0$ $<$ $r_{2}$ $<$ $r_{0},$ such
that

\begin{equation}
|N(q)-N(p)|\leq \frac{1}{2}  \label{3.a}
\end{equation}%
\noindent for all $q$ $\in $ $B_{r_{2}}(p)$. Let $\Gamma _{1},$ $\Gamma _{2}$
are two seed curves passing through $p$. For $j$ $=$ $1,$ $2,$ let $\Gamma
_{j}$ $=$ $\Gamma _{j}^{+}$ $\cup $ $\Gamma _{j}^{-},$ $\Gamma _{j}^{+}$ $%
\cap $ $\Gamma _{j}^{-}$ $=$ $\{p\}$ where $\Gamma _{j}^{+}$ ($\Gamma
_{j}^{-},$ resp.) is the part of $\Gamma _{j}$ emanating from $p$ along the $%
+N$ ($-N,$ resp.) direction. Then for any $r,$ $0$ $\leq $ $r$ $<$ $r_{2},$
there exist $p_{j,r}^{\pm },$ $j$ $=$ $1,$ $2$ such that $\partial B_{r}(p)$ 
$\cap $ $\Gamma _{j}^{\pm }$ $=$ $\{p_{j,r}^{\pm }\}$ where $\partial
B_{0}(p)$ $\equiv $ $\{p\}.$ Suppose (a) fails to hold. Then without loss of
generality, we may assume there exists $r_{4},$ $0$ $<$ $r_{4}$ $<$ $r_{2}$
such that $p_{1,r_{4}}^{+}$ $\neq $ $p_{2,r_{4}}^{+}$ and there exists a
unique $r_{3}$ depending on $r_{4}$ only such that $0$ $\leq $ $r_{3}$ $<$ $%
r_{4}$ and $p_{1,r_{3}}^{+}$ $=$ $p_{2,r_{3}}^{+}$ (if $r_{3}$ $=$ $0,$ $%
p_{1,r_{3}}^{+}$ $=$ $p_{2,r_{3}}^{+}$ $=$ $p),$ $p_{1,r}^{+}$ $\neq $ $%
p_{2,r}^{+}$ for $r_{3}$ $<$ $r$ $<$ $r_{4}$ (see Figure 3.1)$.$ 
\begin{figure}[b]
\centering              
\begin{pspicture}(0,-1.5)(11,4.5) 
\psdots[dotstyle=*,dotscale=1 1](0,0)
\psarc[linewidth=1.0pt]{-}(0 ,4){4}{270}{345}
\psarc[linewidth=1.0pt]{-}(0 ,6){6}{270}{320}
\psarc[linewidth=1.0pt]{-}(0 ,0){4.5}{10}{45}

\psdots[dotstyle=*,dotscale=1 1](6,0)
\psdots[dotstyle=*,dotscale=1 1](8,0.53)
\uput{5pt}[-90]{0}(0,0){$0=r_3$}
\uput{5pt}[0]{0}(4.2,0.4){$r_4$}
\uput{0pt}[0]{0}(5.4,2){\large{or}}

\psarc[linewidth=1.0pt]{-}(6 ,4){4}{270}{345}
\psarc[linewidth=1.0pt]{-}(5 ,5.72){6}{-60}{-30}
\psarc[linewidth=1.0pt]{-}(6 ,0){4.5}{10}{45}
\pscurve[linestyle=dashed, dash=2pt 3pt, linewidth=1.0pt]{-}(6,0)(7.5,0.15)(8,0.52)
\uput{5pt}[-90]{0}(6,0){$0$}
\uput{5pt}[-90]{0}(8,0.4){$r_3$}
\uput{5pt}[0]{0}(10.2,0.4){$r_4$}
\uput{0pt}[0]{0}(4.7,-1.5){Figure 3.1}
\end{pspicture}
\end{figure}

Let $\ell _{r}$ denote the shorter arc of $\partial B_{r}(p)$ connecting $%
p_{1,r}^{+}$ and $p_{2,r}^{+}.$ For perhaps smaller $r_{3},$ $r_{4}$ we have%
\begin{equation}
N(p)\cdot \partial _{r}(q)\geq \frac{3}{4}  \label{3.b}
\end{equation}%
\noindent where $q$ $\in $ $\ell _{r}$ and $\partial _{r}(q)$ $=$ $\frac{q-p%
}{|q-p|}$ is the unit outer normal to $\ell _{r}$ for $r_{3}$ $<$ $r$ $<$ $%
r_{4}.$ It then follows from (\ref{3.a}) and (\ref{3.b}) that 
\begin{equation}
\frac{1}{4}\leq N\cdot \partial _{r}\leq 1\text{ \ on }\ell _{r}  \label{3.1}
\end{equation}

\noindent for $r_{3}$ $<$ $r$ $<$ $r_{4}.$ Let $\Omega _{r}$ be the domain
surrounded by $\Gamma _{i}^{+}$ $\cap $ $(B_{r}(p)\backslash \bar{B}%
_{r_{3}}(p)),$ $i$ $=$ $1,$ $2,$ and $\ell _{r}$ with vertices $p_{3}$ ($=$ $%
p_{1,r_{3}}^{+}$ $=$ $p_{2,r_{3}}^{+}),$ $p_{2,r}^{+},$ and $p_{1,r}^{+},$
where $\bar{B}_{r_{3}}(p)$ $\equiv $ $\{p\}$ if $r_{3}$ $=$ $0.$ Let

\begin{equation}
h(r)\equiv \doint_{\partial \Omega _{r}}N\cdot \nu .  \label{3.2}
\end{equation}%
\noindent Observe that $N\cdot \nu $ $=$ $0$ along $\Gamma _{1}$ and $\Gamma
_{2}$ and $\nu $ $=$ $\partial _{r}$ on $\ell _{r}.$ It follows from (\ref%
{3.1}) and (\ref{3.2}) that

\begin{equation}
\frac{1}{4}|\ell _{r}|\leq h(r)\leq |\ell _{r}|  \label{3.3}
\end{equation}%
\noindent for $r_{3}$ $<$ $r$ $<$ $r_{4}$ where $|\ell _{r}|$ denotes the
arc length of $\ell _{r}.$ On the other hand, we compute from (\ref{2.3})
that%
\begin{equation*}
h(r)\equiv \doint_{\partial \Omega _{r}}N\cdot \nu =\int_{\Omega _{r}}H
\end{equation*}

\noindent and hence by (\ref{3.3}) we have%
\begin{equation}
h^{\prime }(r)=\int_{\ell _{r}}H\leq H_{M}(r)|\ell _{r}|\leq 4H_{M}(r)h(r)
\label{3.4}
\end{equation}

\noindent for $r_{3}$ $<$ $r$ $<$ $r_{4}.$ We can therefore have%
\begin{eqnarray}
&&\frac{d}{dr}[h(r)e^{-4\int_{r_{3}}^{r}H_{M}(r)dr}]  \label{3.5} \\
&=&[h^{\prime }(r)-4H_{M}(r)h(r)]e^{-4\int_{r_{3}}^{r}H_{M}(r)dr}\leq 0 
\notag
\end{eqnarray}

\noindent by (\ref{3.4}). It follows from (\ref{3.5}) that 
\begin{equation}
h(r)e^{-4\int_{r_{3}}^{r}H_{M}(r)dr}\leq h(r_{3})=0.  \label{3.6}
\end{equation}

\noindent By assumption (\ref{1.7}) we have $\int_{0}^{r}H_{M}(r)dr$ $<$ $%
\infty $ and hence $h(r)$ $=$ $0$ by (\ref{3.6}) for $r_{3}$ $\leq $ $r$ $%
\leq $ $r_{4}.$ From (\ref{3.3}) $|\ell _{r}|$ $=$ $0$ for $r_{3}$ $\leq $ $%
r $ $\leq $ $r_{4},$ which implies that $p_{1,r}^{+}$ $=$ $p_{2,r}^{+}$ for $%
r_{3}$ $<$ $r$ $<$ $r_{4}.$ We have reached a contradiction and hence proved
(a).

Next we will prove (b). Suppose $\Gamma _{1}$ and $\Gamma _{2}$ are two
characteristic curves passing through $p.$ Without loss of generality we may
assume that locally $\Gamma _{1}$ and $\Gamma _{2}$ are graphs $(x,$ $%
y_{1}(x))$ and $(x,$ $y_{2}(x)),$ $|x|$ $\leq $ $x_{1}$ for some positive
constant $x_{1},$ respectively and $p$ $=$ $(0,$ $y_{1}(0))$ $=$ $(0,$ $%
y_{2}(0)),$ $y_{1}^{\prime }(0)$ $=$ $y_{2}^{\prime }(0)$ $=$ $0.$ To prove
(b) we need to show that $y_{1}(x)$ $=$ $y_{2}(x)$ on $[-x_{0}^{\prime },$ $%
x_{0}^{\prime }]$ for some small positive number $x_{0}^{\prime }$ $\leq $ $%
x_{1}.$ We will only show that $y_{1}(x)$ $=$ $y_{2}(x)$ on $[0,$ $x_{0}]$
for some small positive number $x_{0}$ $\leq $ $x_{1}$(see below) since a
similar argument works for the interval of nonpositive numbers. We take $%
x_{0}$ $=$ $\frac{1}{2H_{1}}$ for a large constant $H_{1}$ such that $|H|$ $%
\leq $ $H_{1}$ on $\Gamma _{1}$ and $\Gamma _{2}$ for $x$ $\in $ $[0,$ $%
x_{0}].$

By Theorem A the curvature of $\Gamma _{j},$ $j$ $=$ $1,$ $2,$ equals $-H.$
Namely, we have (\ref{2.7}) and hence (say, for $x$ $\in $ $[0,$ $\frac{1}{%
2H_{1}}])$%
\begin{equation}
\frac{y_{j}^{\prime }(x)}{\sqrt{1+(y_{j}^{\prime }(x))^{2}}}%
=-\int_{0}^{x}H(x,y_{j}(x))dx  \label{3.7}
\end{equation}

\noindent for $j$ $=$ $1,$ $2.$ It follows that%
\begin{equation}
|\frac{y_{j}^{\prime }(x)}{\sqrt{1+(y_{j}^{\prime }(x))^{2}}}|\text{ }\leq 
\text{ }\frac{1}{2}  \label{3.8a}
\end{equation}%
\noindent and hence%
\begin{equation}
|y_{j}^{\prime }(x)|\text{ }\leq \text{ }\frac{1}{\sqrt{3}}  \label{3.9a}
\end{equation}%
\noindent for $x$ $\in $ $[0,$ $\frac{1}{2H_{1}}],$ $j$ $=$ $1,$ $2.$ We
want to prove that $y_{2}(x)$ $=$ $y_{1}(x)$ for all $x$ $\in $ $[0,$ $\frac{%
1}{2H_{1}}].$ If not, we may assume 
\begin{equation}
y_{2}(x)>y_{1}(x)  \label{3.9b}
\end{equation}

\noindent\ for all $x$ $\in $ $(0,$ $\frac{1}{2H_{1}}]$ (otherwise we can
find an interval $[a,b]$ $\subset $ $[0,$ $\frac{1}{2H_{1}}],$ $a$ $<$ $b,$
such that $y_{1}(a)$ $=$ $y_{2}(a),$ $y_{1}(b)$ $=$ $y_{2}(b),$ $%
y_{1}^{\prime }(a)$ $=$ $y_{2}^{\prime }(a),$ and $y_{2}(x)>y_{1}(x)$ (or $%
y_{2}(x)<y_{1}(x))$ for $x$ $\in $ $(a,b).$ Applying a similar reasoning
below to $[a,b]$ instead of $[0,$ $\frac{1}{2H_{1}}],$ we will reach $%
y_{2}(x)$ $=$ $y_{1}(x)$ for all $x$ $\in $ $[a,b],$ a contradiction). From (%
\ref{3.7}) we compute%
\begin{eqnarray}
&&|\frac{y_{2}^{\prime }(x)}{\sqrt{1+(y_{2}^{\prime }(x))^{2}}}-\frac{%
y_{1}^{\prime }(x)}{\sqrt{1+(y_{1}^{\prime }(x))^{2}}}|  \label{3.10a} \\
&\leq &\int_{0}^{x}|H(x,y_{2}(x))-H(x,y_{1}(x))|dx\leq
C_{1}\int_{0}^{x}(y_{2}(x)-y_{1}(x))dx  \notag
\end{eqnarray}

\noindent where $C_{1}$ is the Lipschitzian constant of $H.$ Let%
\begin{equation*}
h(x)\equiv \int_{0}^{x}(y_{2}(x)-y_{1}(x))dx.
\end{equation*}

\noindent It follows that $h^{\prime }(x)$ $=$ $y_{2}(x)-y_{1}(x)$ and $%
h^{\prime \prime }(x)$ $=$ $y_{2}^{\prime }(x)$ $-$\ $y_{1}^{\prime }(x).$
On the other hand, observe that $f^{\prime }(t)$ $=$ $(\frac{1}{1+t^{2}}%
)^{3/2}$ for $f(t)$ $=$ $\frac{t}{(1+t^{2})^{1/2}}.$ By the mean-value
theorem and (\ref{3.9a}), we have%
\begin{eqnarray}
&&|\frac{y_{2}^{\prime }(x)}{\sqrt{1+(y_{2}^{\prime }(x))^{2}}}-\frac{%
y_{1}^{\prime }(x)}{\sqrt{1+(y_{1}^{\prime }(x))^{2}}}|  \label{3.11a} \\
&\geq &(\frac{3}{4})^{3/2}|y_{2}^{\prime }(x)-y_{1}^{\prime }(x)|.  \notag
\end{eqnarray}

From (\ref{3.11a}) and (\ref{3.10a}) we obtain the following differential
inequality for $h:$ $(C_{2}$ $=$ $(\frac{4}{3})^{3/2}C_{1})$%
\begin{equation}
h^{\prime \prime }(x)\leq C_{2}h(x).  \label{3.12a}
\end{equation}

\noindent Multiplying (\ref{3.12a}) by $h^{\prime }(x)$ ($>0$ by (\ref{3.9b}%
)) and integrating from $0$ to $x$ $\in $ $(0,$ $\frac{1}{2H_{1}}]$, we get%
\begin{equation}
h^{\prime }(x)\leq \sqrt{C_{2}}h(x)  \label{3.13a}
\end{equation}

\noindent in view of $h(x)$ $>$ $0,$ $h^{\prime }(x)$ $>$ $0,$ and $%
h^{\prime }(0)$ $=$ $h(0)$ $=$ $0.$ Writing (\ref{3.13a}) as $(\log
h)^{\prime }(x)$ $\leq $ $\sqrt{C_{2}}$ and then integrating from $%
\varepsilon $ ($>$ $0)$ to $x,$ we obtain%
\begin{equation}
h(x)\leq h(\varepsilon )e^{\sqrt{C_{2}}(x-\varepsilon )}.  \label{3.14a}
\end{equation}

\noindent Letting $\varepsilon \rightarrow 0$ in (\ref{3.14a}) gives $h(x)$ $%
\equiv $ $0$ on $[0,$ $\frac{1}{2H_{1}}],$ a contradiction. We have proved
(b).

\textbf{\bigskip }%
\endproof%
\endproof%

We remark that we can give an alternative proof of part (b) of Theorem B by
applying Picard-Lindel\"{o}f's Theorem (\cite{Ha}, Theorem 1.1) to (\ref{3.7}%
). Recall (see Section 1) that for $u$ $\in $ $C^{1}$ and $\vec{F}$ $\in $ $%
C^{1},$ let $D$ $\equiv $ $|\nabla u$ $+$ $\vec{F}|$, and if $D$ $\neq $ $0,$
we let $N$ $\equiv $ $N_{\vec{F}}^{u}$ $\equiv $ $\frac{\nabla u+\vec{F}}{%
|\nabla u+\vec{F}|}.$ Hence $N^{\perp }$ $\equiv $ $N_{\vec{F}}^{u,\perp }$ $%
=$ $D^{-1}$ $(u_{y}+F_{2},$ $-u_{x}-F_{1})$ where we write $\vec{F}$ $=$ $%
(F_{1},$ $F_{2}).$ Recall the definition of $rot\vec{F}$ as follows:%
\begin{equation*}
rot\vec{F}=(F_{2})_{x}-(F_{1})_{y}.
\end{equation*}

\noindent It is easy to see that $div$ $(DN^{\perp })$ $=$ $rot\vec{F}$ if $%
u $ $\in $ $C^{2}.$ Note that $rot\vec{F}$ $=$ $2$ for $\vec{F}$ $=$ $(-y,$ $%
x).$ For $u$ $\in $ $C^{1},$ we have the following result.

\bigskip

\textbf{Lemma 3.1}. \textit{Let }$u:\Omega \subset R^{2}\rightarrow R$%
\textit{\ be a }$C^{1}$\textit{\ smooth function such that }$D$\textit{\ }$%
\neq $\textit{\ }$0$\textit{\ on }$\Omega $\textit{\ (i.e., }$\Omega $%
\textit{\ is a nonsingular domain) with }$\vec{F}$ $\in $ $C^{1}(\Omega )$%
\textit{. Let }$\Omega ^{\prime \prime }$\textit{\ }$\subset \subset $%
\textit{\ }$\Omega $\textit{\ be a bounded Lipschitzian domain. Then for }$%
\varphi $\textit{\ }$\in $\textit{\ }$C^{1}(\Omega )$\textit{\ there holds}%
\begin{equation}
\oint_{\partial \Omega ^{\prime \prime }}\varphi DN^{\perp }\cdot \nu
=\int_{\Omega ^{\prime \prime }}\{(\nabla \varphi )\cdot (DN^{\perp
})+\varphi rot\vec{F}\}.  \label{4.7.0}
\end{equation}

\bigskip

\proof
Write $DN^{\perp }$ $=$ $(u_{y},$ $-u_{x})$ $+$ $(F_{2},$ $-F_{1}).$ Let $%
v_{\varepsilon }$ denote a mollifier of $v.$ Observe that $%
(u_{y})_{\varepsilon }$ $=$ $(u_{\varepsilon })_{y},$ $(u_{x})_{\varepsilon
} $ $=$ $(u_{\varepsilon })_{x}.$ It follows that%
\begin{eqnarray}
div((u_{y})_{\varepsilon },(-u_{x})_{\varepsilon }) &=&div((u_{\varepsilon
})_{y},-(u_{\varepsilon })_{x})  \label{4.7.1} \\
&=&(u_{\varepsilon })_{yx}-(u_{\varepsilon })_{xy}=0.  \notag
\end{eqnarray}

Now using the divergence theorem, we compute%
\begin{eqnarray}
&&\oint_{\partial \Omega ^{\prime \prime }}\varphi \lbrack
((u_{y})_{\varepsilon },(-u_{x})_{\varepsilon })+(F_{2},-F_{1})]\cdot \nu
\label{4.7.2} \\
&=&\int_{\Omega ^{\prime \prime }}(\nabla \varphi )\cdot \lbrack
((u_{y})_{\varepsilon },(-u_{x})_{\varepsilon })+(F_{2},-F_{1})]+\varphi
div[((u_{y})_{\varepsilon },(-u_{x})_{\varepsilon })+(F_{2},-F_{1})]  \notag
\\
&=&\int_{\Omega ^{\prime \prime }}(\nabla \varphi )\cdot \lbrack
((u_{y})_{\varepsilon },(-u_{x})_{\varepsilon })+(F_{2},-F_{1})]+\varphi rot%
\vec{F}  \notag
\end{eqnarray}

\noindent by (\ref{4.7.1}) and noting that $div(F_{2},-F_{1})$ $=$ $rot\vec{F%
}.$ Taking the limit $\varepsilon $ $\rightarrow $ $0$ in (\ref{4.7.2})
gives (\ref{4.7.0}).

\bigskip 
\endproof%

\proof
(of Theorem B$^{\prime }$) Since $N^{\perp }$ $\in $ $C^{0}(\Omega ),$ we
can choose $r_{2}$ $>$ $0$ such that $B_{r_{2}}(p)$ $\subset \subset $ $%
\Omega $ and 
\begin{equation}
|N^{\perp }(q)-N^{\perp }(p)|\leq \frac{1}{2}  \label{3.14}
\end{equation}

\noindent for all $q$ $\in $ $B_{r_{2}}(p)$. Let $\Gamma _{1},$ $\Gamma _{2}$
be two characteristic curves passing through $p$. For $j$ $=$ $1,$ $2,$ let $%
\Gamma _{j}$ $=$ $\Gamma _{j}^{+}$ $\cup $ $\Gamma _{j}^{-},$ $\Gamma
_{j}^{+}$ $\cap $ $\Gamma _{j}^{-}$ $=$ $\{p\}$ where $\Gamma _{j}^{+}$ ($%
\Gamma _{j}^{-},$ resp.) is the part of $\Gamma _{j}$ emanating from $p$
along the $+N^{\perp }$ ($-N^{\perp },$ resp.) direction.Then for any $r,$ $%
0 $ $\leq $ $r$ $<$ $r_{2},$ there exist $p_{j,r}^{\pm },$ $j$ $=$ $1,$ $2$
such that $\partial B_{r}(p)$ $\cap $ $\Gamma _{j}^{\pm }$ $=$ $%
\{p_{j,r}^{\pm }\}.$ Suppose the conclusion is false. Then without loss of
generality, we may assume there exists $r_{4},$ $0$ $<$ $r_{4}$ $<$ $r_{2}$
such that $p_{1,r_{4}}^{+}$ $\neq $ $p_{2,r_{4}}^{+}$ and there exists a
unique $r_{3}$ depending on $r_{4}$ only such that $0$ $\leq $ $r_{3}$ $<$ $%
r_{4}$ and $p_{1,r_{3}}^{+}$ $=$ $p_{2,r_{3}}^{+}$ (if $r_{3}$ $=$ $0,$ $%
p_{1,r_{3}}^{+}$ $=$ $p_{2,r_{3}}^{+}$ $=$ $p),$ $p_{1,r}^{+}$ $\neq $ $%
p_{2,r}^{+}$ for $r_{3}$ $<$ $r$ $<$ $r_{4}$. Let $\ell _{r}$ denote the
shortest arc of $\partial B_{r}(p)$ connecting $p_{1,r}^{+}$ and $%
p_{2,r}^{+}.$ For perhaps smaller $r_{3},$ $r_{4}$ we have%
\begin{equation}
N^{\perp }(p)\cdot \partial _{r}(q)\geq \frac{3}{4}  \label{3.14.1}
\end{equation}%
\noindent where $q$ $\in $ $\ell _{r}$ and $\partial _{r}(q)$ $=$ $\frac{q-p%
}{|q-p|}$ is the unit outer normal to $\ell _{r}$ for $r_{3}$ $<$ $r$ $<$ $%
r_{4}.$ It then follows from (\ref{3.14}) and (\ref{3.14.1}) that%
\begin{equation}
\frac{1}{4}\leq N^{\perp }\cdot \partial _{r}\leq 1\text{ \ on }\ell _{r}
\label{3.14.2}
\end{equation}

\noindent for $r_{3}$ $<$ $r$ $<$ $r_{4}.$ Let $\Omega _{r}$ be the domain
surrounded by $\Gamma _{i}^{+}$ $\cap $ $(B_{r}(p)\backslash \bar{B}%
_{r_{3}}(p)),$ $i$ $=$ $1,$ $2,$ and $\ell _{r}$ with vertices $p_{3}$ ($=$ $%
p_{1,r_{3}}^{+}$ $=$ $p_{2,r_{3}}^{+}),$ $p_{2,r}^{+},$ and $p_{1,r}^{+}.$
Let%
\begin{equation}
h(r)\equiv \doint_{\partial \Omega _{r}}DN^{\perp }\cdot \nu .  \label{3.15}
\end{equation}

\noindent Observe that $N^{\perp }\cdot \nu $ $=$ $0$ along $\Gamma _{1}$
and $\Gamma _{2}$ and $\nu $ $=$ $\partial _{r}$ on $\ell _{r}.$ It follows
from (\ref{3.14.2}) and (\ref{3.15}) that%
\begin{equation}
\frac{1}{4}C_{1}|\ell _{r}|\leq h(r)\leq C_{2}|\ell _{r}|  \label{3.16}
\end{equation}

\noindent where $|\ell _{r}|$ denotes the arc length of $\ell _{r},$ $C_{1}$ 
$\equiv $ $\min_{\bar{B}_{r_{4}}(p)}D$ $>$ $0$ since $D$ $>$ $0$ on the
nonsingular domain $\Omega ,$ and $C_{2}$ $\equiv $ $\max_{\bar{B}%
_{r_{4}}(p)}D.$ On the other hand, we compute from (\ref{3.15}) and (\ref%
{4.7.0}) with $\varphi $ $\equiv $ $1$ that%
\begin{equation*}
h(r)\equiv \doint_{\partial \Omega _{r}}DN^{\perp }\cdot \nu =\int_{\Omega
_{r}}rot\vec{F}
\end{equation*}

\noindent and hence 
\begin{equation}
h^{\prime }(r)=\int_{\ell _{r}}rot\vec{F}\leq C_{3}|\ell _{r}|\leq \frac{%
4C_{3}}{C_{1}}h(r)  \label{3.17}
\end{equation}

\noindent by (\ref{3.16}) for $r_{3}$ $<$ $r$ $<$ $r_{4},$ where $|rot\vec{F}%
|$ $\leq $ $C_{3}$ $\equiv $ $\max_{\bar{B}_{r_{4}}(p)}|rot\vec{F}|.$ We can
therefore have%
\begin{eqnarray}
&&\frac{d}{dr}[h(r)e^{-\frac{4C_{3}r}{C_{1}}}]  \label{3.18} \\
&=&[h^{\prime }(r)-\frac{4C_{3}}{C_{1}}h(r)]e^{-\frac{4C_{3}r}{C_{1}}}\leq 0
\notag
\end{eqnarray}

\noindent by (\ref{3.17}). It follows from (\ref{3.18}) that%
\begin{equation*}
h(r)e^{-\frac{4C_{3}r}{C_{1}}}\leq h(r_{3})e^{-\frac{4C_{3}r_{3}}{C_{1}}}=0.
\end{equation*}

\noindent Therefore $h(r)$ $=$ $0$ for $r_{3}$ $<$ $r$ $<$ $r_{4}$ and then
from (\ref{3.16}) we have $|\ell _{r}|$ $=$ $0$ for $r_{3}$ $<$ $r$ $<$ $%
r_{4}$ which implies that $p_{1,r}^{+}$ $=$ $p_{2,r}^{+}$ for $r_{3}$ $<$ $r$
$<$ $r_{4}$, a contradiction$.$

\bigskip 
\endproof%

Note that on the boundary of any Lipschitzian domain $\Omega ^{\prime }$ $%
\subset \subset $ $\Omega $ ($\Omega $ being nonsingular), $N^{\perp }$ $%
\cdot $ $\nu $ $=$ $N$ $\cdot $ $(\frac{dx}{d\sigma },\frac{dy}{d\sigma })$
where $\sigma $ is the unit-speed parameter and the unit tangent $(\frac{dx}{%
d\sigma },\frac{dy}{d\sigma })$ is the rotation of the unit outer normal $%
\nu $ by $\frac{\pi }{2}$ degrees. It follows that for $\vec{F}$ $=$ $(-y,$ $%
x),$ $DN^{\perp }$ $\cdot $ $\nu $ $d\sigma $ $=$ $DN$ $\cdot $ $(dx,dy)$ $=$
$(u_{x}-y)$ $dx$ $+$ $(u_{y}+x)$ $dy$ $=$ $du$ $+$ $xdy$ $-$ $ydx$ which is
the standard contact form of $R^{3}$, restricted to the surface $(x,$ $y,$ $%
u(x,y)).$

\bigskip

\textbf{Example 3.2.} We will define a family of curves to be the
characteristic curves (for $N^{\perp }$ being its unit tangent vector
field). Let $p_{0}$ $=$ $(x_{0},$ $y_{0})$ in the $xy$-plane. Case 1: For $%
y_{0}$ $-$ $x_{0}^{4}$ $\geq $ $0$, we take $y$ $=$ $x^{4}$ $+$ $(y_{0}$ $-$ 
$x_{0}^{4})$ to be the characteristic curve passing through $p_{0}$. Case 2:
For $y_{0}$ $-$ $x_{0}^{4}$ $<$ $0$ and $y_{0}$ $>$ $0$ (or equivalently, $0$
$<$ $\frac{y_{0}}{x_{0}^{4}}$ $<$ $1$ and $x_{0}$ $\neq $ $0),$ we take $y$ $%
=$ $\frac{y_{0}}{x_{0}^{4}}x^{4}$ to be the characteristic curve passing
through $p_{0}.$ Case 3: For $y_{0}$ $\leq $ $0,$ we take $y$ $=$ $y_{0}$ to
be the characteristic curve passing through $p_{0}.$ Note that there are
infinite number of the above-mentioned (characteristic) curves passing
through the origin $(0,$ $0)$ (see Figure 3.2 below)$.$

\bigskip

\begin{figure}[tph]
\centering     
\begin{pspicture}(-3,-2)(3,5) 
\psset{algebraic=true, xunit=1.5,yunit=0.5, linewidth=0.5pt} 
\psaxes[ticks=none, labels=none]{->}(0,0)(-2,-3.5)(2,9) 
\psplot[algebraic=true,linewidth=1.0pt]{-1.5}{1.5}{x^4+3} 
\psplot[algebraic=true,linewidth=1.0pt]{-1.5}{1.5}{x^4+1.5}
\psplot[algebraic=true,linewidth=1.0pt]{-1.5}{1.5}{x^4}
\psplot[algebraic=true,linewidth=1.0pt]{-1.5}{1.5}{0.5*x^4}
\psplot[algebraic=true,linewidth=1.0pt]{-1.5}{1.5}{0.2*x^4}
\psplot[algebraic=true,linewidth=1.0pt]{-1.5}{1.5}{-1} 
\psplot[algebraic=true,linewidth=1.0pt]{-1.5}{1.5}{-2}
\psplot[algebraic=true,linewidth=1.0pt]{-1.5}{1.5}{-3} 

\uput{3pt}[90]{0}(0,9){$y$}
\uput{3pt}[0]{0}(2,0){$x$}
\uput{3pt}[0]{0}(1.5,7.5){$y=x^4+C_1, C_1>0$}
\uput{3pt}[0]{0}(1.5,4.5){$y=x^4$}
\uput{3pt}[0]{0}(1.5,1.5){$y=ax^4, 0<a<1$}
\uput{3pt}[0]{0}(1.5,-2){$y=C_2, C_2<0$}

\uput*{4pt}[-90]{0}(0,-4){Figure 3.2} 
\end{pspicture}
\end{figure}

\bigskip

We compute the curvature of this family of curves as follows. For case 1,
the curvature at $p_{0}$ equals $12x_{0}^{2}(1+16x_{0}^{6})^{-3/2}$. For
case 2, the curvature at $p_{0}$ equals $12(\frac{y_{0}}{x_{0}^{4}}%
)x_{0}^{2}[1+16(\frac{y_{0}}{x_{0}^{4}})^{2}x_{0}^{6}]^{-3/2}.$ For case 3,
the curvature at $p_{0}$ vanishes. According to Theorem A the curvature of a
characteristic curve is $-H$ (with this value (\ref{1.6}) holds since $N$ $=$
$(\cos \theta ,$ $\sin \theta )$ is $C^{1}$ smooth away from $(0,$ $0)).$ So
we can easily verify that $H$ $\in $ $C^{0}.$ On the other hand, we observe
that 
\begin{equation*}
\frac{H(x,y)-H(x,0)}{y}=-\frac{12x^{2}(1+16x^{6})^{-3/2}-0}{x^{4}}%
\rightarrow -\infty
\end{equation*}

\noindent as $x\rightarrow 0$ for $y$ $=$ $x^{4}.$ Therefore $H$ $\notin $ $%
C^{0,1}$ in a neighborhood of $(0,$ $0).$ Altogether this is a
counterexample to Theorem B (b) if $H$ is only continuous, but not
Lipschitzian.

\bigskip

\section{Characteristic coordinates}

In this section we will introduce a system of characteristic coordinates
(see Definition 1.2) for later use.

\bigskip

\textbf{Theorem 4.1}. \textit{Let }$\Omega $\textit{\ be a domain of }$R^{2}$%
\textit{\ and }$H\in C^{0}(\Omega ).$\textit{\ Let }$\theta \in C^{0}(\Omega
)$\textit{\ satisfy (\ref{1.6}), i.e.}%
\begin{equation}
\int_{\Omega }(\cos \theta ,\sin \theta )\cdot \nabla \varphi +\int_{\Omega
}H\varphi =0  \label{4.1}
\end{equation}

\noindent \textit{for all }$\varphi \in C_{0}^{\infty }(\Omega ).$\textit{\
Then given a point }$p_{0}$\textit{\ }$\in $\textit{\ }$\Omega ,$\textit{\
there exist a small neighborhood }$\Omega ^{\prime }$\textit{\ of }$p_{0}$%
\textit{\ and a function }$s$\textit{\ }$\in $\textit{\ }$C^{1}(\Omega
^{\prime })$\textit{\ such that }$\nabla s$\textit{\ }$=$\textit{\ }$%
fN^{\perp }$\textit{\ for some positive function }$f$\textit{\ }$\in $%
\textit{\ }$C^{0}(\Omega ^{\prime })$\textit{\ and the curves defined by }$s$%
\textit{\ }$=$\textit{\ }$c,$ \textit{constants, are seed curves}$.$\textit{%
\ Moreover, }$Nf$ exists and is continuous. In fact, $f$\textit{\ satisfies
the following equation:}%
\begin{equation}
Nf+fH=0.  \label{4.1.1}
\end{equation}

\bigskip

\proof
Without loss of generality we may assume $p_{0}$ $=$ $(0,0),$ the origin,
and $\theta (p_{0})$ $=$ $\frac{\pi }{2}.$ That is, $N$ $\equiv $ $(\cos
\theta ,$ $\sin \theta )$ $=$ $(0,$ $1)$ at $p_{0}.$ Let $\Upsilon $ denote
the $x$-axis$,$ i.e., $\Upsilon $ $=$ $\{$ $(x,y)$ $\in $ $\Omega $ $|$ $y=0$
$\}.$ Since $\theta $ is continuous, we can find a small ball $%
B_{r_{1}}(p_{0})$ $\subset \subset $ $\Omega $ of radius $r_{1}$ $>$ $0$
such that for any point $q$ $\in $ $B_{r_{1}}(p_{0}),$ $|\theta (q)$ $-$ $%
\frac{\pi }{2}|$ $<<1$ and there exists a seed curve passing through $q$ and
intersecting $\Upsilon $ at $p.$ Now we define $s$ $:$ $B_{r_{1}}(p_{0})$ $%
\rightarrow $ $R$ by%
\begin{equation}
s(q)=s(p)=x\text{ if }p=(x,0).  \label{4.2}
\end{equation}

\noindent Then $s$ is well defined by the uniqueness of seed curves
according to Theorem B (a). We can find smaller positive numbers $r_{3}$ $<\
r_{2}$ $<$ $r_{1}$ such that $B_{r_{3}}(p_{0})$ $\subset \subset $ $%
(-r_{2},r_{2})$ $\times $ $(-r_{2},r_{2})$ $\subset \subset $ $%
B_{r_{1}}(p_{0})$ and for $c$ $\in $ $R,$ if $\{s=c\}$ $\cap $ $%
B_{r_{3}}(p_{0})$ $\neq $ $\phi ,$ then $\{s=c\}$ $\cap $ $(-r_{2},r_{2})$ $%
\times $ $(-r_{2},r_{2})$ is a graph, denoted by $(x^{c}(y),$ $y),$ for $%
-r_{2}$ $<$ $y$ $<$ $r_{2}$ (see Figure 4.1)$.$

\begin{figure}[h]
\centering              
\begin{pspicture}(-1,-1.5)(6.5,5.5) 
\psaxes[linewidth=0.5pt,labels=none,ticks=none,Ox=0,Oy=0,Dx=1,Dy=1]{->}(0,0)(-0.5,-0.5)(6,5)
\pscurve[linewidth=1.0pt]{-}(2,0)(1.85,1)(2,2)(2.2,3)(2.15,4)
\psdots[dotstyle=*,dotscale=1 1](2.15,4)
\psline[linewidth=1.0pt]{->}(1.88,0.6)(1.62,2)
\psdots[dotstyle=*,dotscale=0.5 0.5](1.88,0.6)
\uput{0pt}[0]{0}(1.3,1.6){\tiny{$N$}}

\pscurve[linewidth=1.0pt]{-}(4,0)(3.92,1)(4,2)(4.2,3)(4.15,4)
\psdots[dotstyle=*,dotscale=1 1](4.15,4)
\psline[linewidth=1.0pt]{->}(3.94,0.6)(3.79,2)
\psdots[dotstyle=*,dotscale=0.5 0.5](3.94,0.6)
\uput{0pt}[0]{0}(3.3,1.6){\tiny{$N$}}
\uput{0pt}[0]{0}(2.6,1.8){$\stackrel{\thicksim}{\Omega}$}

\uput{5pt}[90]{0}(2.15,4){\tiny{$(x^{s_1}(y),y$)}}
\uput{5pt}[-90]{0}(2,0){\tiny{$s_1$}}

\uput{5pt}[90]{0}(4.15,4){\tiny{$(x^{s_2}(y),y$)}}
\uput{5pt}[-90]{0}(4,0){\tiny{$s_2$}}

\uput{5pt}[0]{0}(6,0){\normalsize{$x$}}
\uput{5pt}[90]{0}(0,5){\normalsize{$y$}}

\uput{5pt}[225]{0}(0,0){\tiny{$p_0=(0,0)$}}

\psset{linewidth=0.7pt}
\qline(-0.2,4)(5.5,4)
\uput{7pt}[180]{0}(0,4){\tiny{$y$}}

\uput{0pt}[0]{0}(2,-1.5){Figure 4.1}
\end{pspicture}
\end{figure}

It follows that for $c_{1}$ $<$ $c_{2}$%
\begin{eqnarray*}
&&s^{-1}((c_{1},c_{2}))\cap B_{r_{3}}(p_{0}) \\
&=&[\bigcup_{c_{1}<c<c_{2}}\{(x^{c}(y),y)|-r_{2}<y<r_{2}\}]\cap
B_{r_{3}}(p_{0})
\end{eqnarray*}

\noindent is open. So $s$ $\in $ $C^{0}(B_{r_{3}}(p_{0}))$ and $x^{c_{4}}(y)$
$>$ $x^{c_{3}}(y)$ if and only if $c_{4}$ $>$ $c_{3}.$ Next we are going to
prove $s$ $\in $ $C^{1}.$ Given a point $(x_{1},$ $y)$ $\in $ $%
B_{r_{3}}(p_{0}),$ let $s_{1}$ $=$ $s(x_{1},$ $y).$ Then $x_{1}$ $=$ $%
x^{s_{1}}(y)$ by the definition of $x^{c}(y).$ Similarly for $(x_{2},$ $y)$
near $(x_{1},$ $y)$ and $x_{2}$ $>$ $x_{1},$ let $s_{2}$ $=$ $s(x_{2},$ $y).$
Then $x_{2}$ $=$ $x^{s_{2}}(y)$ and $s_{2}$ $>$ $s_{1}.$ We want to compute
the difference quotient%
\begin{equation}
\frac{s(x_{2},y)-s(x_{1},y)}{x_{2}-x_{1}}=\frac{s_{2}-s_{1}}{x_{2}-x_{1}}=%
\frac{s_{2}-s_{1}}{x^{s_{2}}(y)-x^{s_{1}}(y)}.  \label{4.3}
\end{equation}

Now let 
\begin{equation}
A(y)\equiv \int_{x^{s_{1}}(y)}^{x^{s_{2}}(y)}\sin \theta (x,y)dx.
\label{4.6}
\end{equation}%
\noindent We want to know the relation between $A(y)$ and $A(0).$ Without
loss of generality we may assume that $y$ $>$ $0.$ Let $\tilde{\Omega}$ $%
\equiv $ $\{$ $(\varsigma ,$ $\eta )$ $|$ $0$ $<$ $\eta $ $<$ $y,$ $%
x^{s_{1}}(\eta )$ $<$ $\varsigma $ $<$ $x^{s_{2}}(\eta )\}.$ Recall that (%
\ref{2.3}) (which is obtained from (\ref{4.1})) reads 
\begin{equation}
\int_{\partial \tilde{\Omega}}N\cdot \nu =\int_{\tilde{\Omega}}H.
\label{4.4}
\end{equation}

\noindent Observe that $\nu $ $=$ $\pm N^{\perp }$ along the seed curves $%
\{s=s_{j}\},$ $j$ $=$ $1,$ $2$ while $\nu $ $=$ $(0,$ $1)$ on $\{$ $(x,y)$ $%
| $ $x^{s_{1}}(y)$ $<$ $x$ $<$ $x^{s_{2}}(y)\}$ and $\nu $ $=$ $(0,$ $-1)$
on $\{$ $(x,$ $0)$ $|$ $s_{1}$ $<$ $x$ $<$ $s_{2}\}.$ It follows from (\ref%
{4.4}) that (recall that $N$ $=$ $(\cos \theta ,$ $\sin \theta ))$%
\begin{equation}
\int_{x^{s_{1}}(y)}^{x^{s_{2}}(y)}\sin \theta
(x,y)dx+\int_{s_{1}}^{s_{2}}(-\sin \theta
(x,0))dx=\int_{0}^{y}(\int_{x^{s_{1}}(\eta )}^{x^{s_{2}}(\eta )}H(\varsigma
,\eta )d\varsigma )d\eta .  \label{4.5}
\end{equation}

\noindent Then by (\ref{4.5}) we deduce that%
\begin{equation}
A^{\prime }(y)=\int_{x^{s_{1}}(y)}^{x^{s_{2}}(y)}H(x,y)dx.  \label{4.7}
\end{equation}

\noindent By the mean value theorem, there exist $\varsigma _{j}$ $=$ $%
\varsigma _{j}(y,$ $s_{1},$ $s_{2}),$ $j$ $=$ $1,$ $2$ such that $%
x^{s_{1}}(y)$ $<$ $\varsigma _{j}$ $<$ $x^{s_{2}}(y)$ and 
\begin{eqnarray}
A(y) &=&(x^{s_{2}}(y)-x^{s_{1}}(y))\sin \theta (\varsigma _{1},y)
\label{4.8} \\
A^{\prime }(y) &=&(x^{s_{2}}(y)-x^{s_{1}}(y))H(\varsigma _{2},y)  \notag
\end{eqnarray}

\noindent in view of (\ref{4.6}) and (\ref{4.7}). By (\ref{4.8}) we obtain%
\begin{equation}
\frac{d\log A(y)}{dy}=\frac{A^{\prime }(y)}{A(y)}=\frac{H(\varsigma _{2},y)}{%
\sin \theta (\varsigma _{1},y)}  \label{4.9}
\end{equation}

\noindent (noting that $\sin \theta $ is close to $1$ near $p_{0}$ where $%
\theta $ equals $\frac{\pi }{2}$ by assumption). Integrating both sides of (%
\ref{4.9}) gives%
\begin{equation}
\frac{A(y)}{A(0)}=\exp \int_{0}^{y}\frac{H(\varsigma _{2},\eta )}{\sin
\theta (\varsigma _{1},\eta )}d\eta .  \label{4.10}
\end{equation}

\noindent Observe that $A(0)$ $=$ $(s_{2}$ $-$ $s_{1})\sin \theta (\varsigma
_{1},0)$ by (\ref{4.8}). It then follows from (\ref{4.10}) that%
\begin{equation}
\frac{s_{2}-s_{1}}{x^{s_{2}}(y)-x^{s_{1}}(y)}=\frac{\sin \theta (\varsigma
_{1}(y),y)}{\sin \theta (\varsigma _{1}(0),0)}\exp (-\int_{0}^{y}\frac{%
H(\varsigma _{2}(\eta ),\eta )}{\sin \theta (\varsigma _{1}(\eta ),\eta )}%
d\eta ).  \label{4.11}
\end{equation}

\noindent We have omitted the dependence of $s_{1}$ and $s_{2}$ for the
expression of $\varsigma _{1}$ and $\varsigma _{2}$ in (\ref{4.11}).
Combining (\ref{4.3}) and (\ref{4.11}) and taking the limit $x_{2}$ $%
\rightarrow $ $x_{1},$ we finally obtain%
\begin{equation}
\frac{\partial s}{\partial x}(x_{1},y)=\frac{\sin \theta (x_{1},y)}{\sin
\theta (s_{1},0)}\exp (-\int_{0}^{y}\frac{H(x^{s_{1}}(\eta ),\eta )}{\sin
\theta (x^{s_{1}}(\eta ),\eta )}d\eta ).  \label{4.12}
\end{equation}

\noindent Here we have applied Lebesgue's Dominated Convergence Theorem
since $\frac{H(\varsigma _{2}(\eta ),\eta )}{\sin \theta (\varsigma
_{1}(\eta ),\eta )}$ is uniformly bounded and $\varsigma _{j}$ converges to $%
x^{s_{1}}$ pointwise$.$ Let $\tau $ denote the arc length (unit-speed)
parameter of the seed curve $\{s=s_{1}\}.$ By (\ref{4.2}), the definition of 
$s,$ we have $\frac{\partial s}{\partial \tau }$ $=$ $0.$ Note that $Nd\tau $
$=$ $(dx,$ $dy)$ along $\{s$ $=$ $s_{1}\}$ $=$ $\{$ $(x^{s_{1}}(y),$ $y)$ $%
\}.$ It then follows that%
\begin{eqnarray}
\frac{\partial s}{\partial y}(x_{1},y) &=&-\frac{\partial s}{\partial x}%
(x_{1},y)\frac{dx^{s_{1}}(y)}{dy}=-\frac{\partial s}{\partial x}(x_{1},y)%
\frac{\cos \theta (x_{1},y)}{\sin \theta (x_{1},y)}  \label{4.13} \\
&=&f(x_{1},y)(-\cos \theta (x_{1},y))  \notag
\end{eqnarray}

\noindent by (\ref{4.12}), where 
\begin{equation}
f(x_{1},y)\equiv \frac{1}{\sin \theta (s_{1},0)}\exp (-\int_{0}^{y}\frac{%
H(x^{s_{1}}(\eta ),\eta )}{\sin \theta (x^{s_{1}}(\eta ),\eta )}d\eta ).
\label{4.14}
\end{equation}

\noindent Suppose $\tau $ $=$ $0$ at $(s_{1},$ $0)$ and $\tau $ $=$ $l$ at $%
(x_{1},$ $y)$. Recall that $Nd\tau $ $=$ $(dx,$ $dy)$ and hence $\sin \theta
(x^{s_{1}}(\eta ),$ $\eta )$ $d\tau $ $=$ $d\eta .$ So we can rewrite (\ref%
{4.14}) as%
\begin{equation}
f(x_{1},y)=f(s_{1},0)\exp (-\int_{0}^{l}H(\Lambda _{(s_{1},0)}(\tau ))d\tau )
\label{4.15}
\end{equation}

\noindent where $\Lambda _{(s_{1},0)}$ denotes the seed curve $\{s$ $=$ $%
s_{1}\}$ from $(s_{1},$ $0)$ to $(x_{1},$ $y),$ parametrized by $\tau .$
Since $s$ $\in $ $C^{0}$, $\frac{H(x^{s_{1}}(\eta ),\eta )}{\sin \theta
(x^{s_{1}}(\eta ),\eta )}$ is uniformly bounded, and $\frac{H(x^{s_{2}}(\eta
),\eta )}{\sin \theta (x^{s_{2}}(\eta ),\eta )}$ converges to $\frac{%
H(x^{s_{1}}(\eta ),\eta )}{\sin \theta (x^{s_{1}}(\eta ),\eta )}$ pointwise
as $x_{2}$ $\rightarrow $ $x_{1},$ we can apply Lebesgue's Dominated
Convergence Theorem to conclude that $f$ is continuous in $x_{1}$ in view of
(\ref{4.14}). On the other hand, $f$ is continuous along the seed curve in
view of (\ref{4.15}). Together we have $f$ $\in $ $C^{0}$ near $p_{0}.$ From
(\ref{4.13}), (\ref{4.12}), and (\ref{4.14}), we have proved $s$ $\in $ $%
C^{1}$ and

\begin{equation*}
\nabla s=f(\sin \theta ,-\cos \theta )=fN^{\perp }
\end{equation*}

\noindent with $f$ $>$ $0$ and $f$ $\in $ $C^{0}$ near $p_{0}.$ Moreover,
recall that $(Nf)(x_{1},y)$ $\equiv $ $df(\Lambda _{(s_{1},0)}(\tau ))/d\tau 
$ at $\tau $ $=$ $l$ (if exists)$.$ Now (\ref{4.1.1}) easily follows from (%
\ref{4.15}).

\endproof%

Note that we do not need $"u"$ (solution to $divN^{u}$ $=$ $H,$ see (\ref%
{1.2})) to construct $s$ in Theorem 4.1$.$

\bigskip

\textbf{Theorem 4.2}. \textit{Let }$u:\Omega \subset R^{2}\rightarrow R$%
\textit{\ be a }$C^{1}$\textit{\ smooth function and }$\vec{F}$ $\in $ $%
C^{1}(\Omega )$ \textit{such that }$D$\textit{\ }$\equiv $ $|\nabla u$ $+$ $%
\vec{F}|$ $\neq $\textit{\ }$0$\textit{\ on }$\Omega $\textit{\ (i.e., }$%
\Omega $\textit{\ is a nonsingular domain). Suppose we have uniqueness of
characteristic curves passing through a common point. Then given a point }$%
p_{0}$\textit{\ }$\in $\textit{\ }$\Omega ,$\textit{\ there exist a
neighborhood }$\Omega ^{\prime }$\textit{\ of }$p_{0}$\textit{\ and a
function }$t$\textit{\ }$\in $\textit{\ }$C^{1}(\Omega ^{\prime })$\textit{\
such that }$\nabla t$\textit{\ }$=$\textit{\ }$gDN(u)$\textit{\ for some
positive }$g$\textit{\ }$\in $\textit{\ }$C^{0}(\Omega ^{\prime })$\textit{\
and the curves defined by }$t$\textit{\ }$=$\textit{\ }$c,$\textit{\
constants, are characteristic curves}$.$\textit{\ Moreover, }$N^{\perp }g$%
\textit{\ exists and is continuous. In fact, }$g$ \textit{satisfies the
following equation:}%
\begin{equation}
N^{\perp }g+\frac{(rot\vec{F})g}{D}=0.  \label{4.16}
\end{equation}%
\textit{\ }

\bigskip

\proof
The idea similar to that in the proof of Theorem 4.1 works by switching the
role of $N$ and $N^{\perp }.$ Let us sketch a proof. Without loss of
generality we may assume $p_{0}$ $=$ $(0,0),$ the origin, and $\theta
(p_{0}) $ $=$ $\frac{\pi }{2}.$ That is, $N^{\perp }$ $\equiv $ $(\sin
\theta ,$ $-\cos \theta )$ $=$ $(1,0)$ at $p_{0}.$ Let $\Upsilon $ denote
the $y$-axis$. $ Since $\theta $ is continuous, we can find a small ball $%
B_{r_{1}}(p_{0})$ of radius $r_{1}$ $>$ $0$ such that for any point $q$ $\in 
$ $B_{r_{1}}(p_{0})$ there exists a characteristic curve passing through $q$
and intersecting with $\Upsilon $ at $p.$ Now we define $t$ $:$ $%
B_{r_{1}}(p_{0})$ $\rightarrow $ $R$ by%
\begin{equation}
t(q)=t(p)=y\text{ if }p=(0,y).  \label{4.17}
\end{equation}

Then $t$ is well defined by the uniqueness of characteristic curves
according to the assumption. We can find smaller positive numbers $r_{3}$ $%
<\ r_{2}$ $<$ $r_{1}$ such that $B_{r_{3}}(p_{0})$ $\subset \subset $ $%
(-r_{2},r_{2})$ $\times $ $(-r_{2},r_{2})$ $\subset \subset $ $%
B_{r_{1}}(p_{0})$ and for $c$ $\in $ $R,$ if $\{t=c\}$ $\cap $ $%
B_{r_{3}}(p_{0})$ $\neq $ $\phi ,$ then $\{t=c\}$ $\cap $ $(-r_{2},r_{2})$ $%
\times $ $(-r_{2},r_{2})$ is a graph, denoted by $(x,$ $y^{c}(x)),$ for $%
-r_{2}$ $<$ $x$ $<$ $r_{2}.$ By a similar argument as in the proof of
Theorem 4.1 we can prove that $t$ is $C^{0}$ in $B_{r_{3}}(p_{0}).$ Given a
point $(x,$ $y_{1})$ $\in $ $B_{r_{3}}(p_{0}),$ let $t_{1}$ $=$ $t(x,y_{1}).$
Then $y_{1}$ $=$ $y^{t_{1}}(x).$ Similarly for $(x,$ $y_{2})$ near $(x,$ $%
y_{1})$ and $y_{2}$ $>$ $y_{1},$ let $t_{2}$ $=$ $t(x,y_{2}).$ Then $y_{2}$ $%
=$ $y^{t_{2}}(x)$ and $t_{2}$ $>$ $t_{1}.$ We want to compute%
\begin{equation*}
\frac{t(x,y_{2})-t(x,y_{1})}{y_{2}-y_{1}}=\frac{t_{2}-t_{1}}{y_{2}-y_{1}}=%
\frac{t_{2}-t_{1}}{y^{t_{2}}(x)-y^{t_{1}}(x)}
\end{equation*}%
\noindent Now let 
\begin{equation}
B(x)\equiv \int_{y^{t_{1}}(x)}^{y^{t_{2}}(x)}D(x,y)\sin \theta (x,y)dy.
\label{4.20}
\end{equation}%
\noindent We want to know how $B(x)$ is related to $B(0).$ Without loss of
generality we may assume that $x$ $>$ $0.$ Instead of (\ref{4.4}) we have%
\begin{equation}
\oint_{\partial \tilde{\Omega}}DN^{\perp }\cdot \nu =\int_{\tilde{\Omega}}rot%
\vec{F}  \label{4.18}
\end{equation}

\noindent by letting $\varphi $ $=$ $1$ in (\ref{4.7.0}), where $\tilde{%
\Omega}$ $\equiv $ $\{$ $(\varsigma ,$ $\eta )$ $|$ $0$ $<$ $\varsigma $ $<$ 
$x,$ $y^{t_{1}}(\varsigma )$ $<$ $\eta $ $<$ $y^{t_{2}}(\varsigma )\}$.
Observe that $\nu $ $=$ $\pm N$ along the characteristic curves $%
\{t=t_{j}\}, $ $j$ $=$ $1,$ $2$ while $\nu $ $=$ $(1,$ $0)$ on $\{$ $(x,y)$ $%
|$ $y^{t_{1}}(x)$ $<$ $y$ $<$ $y^{t_{2}}(x)\}$ and $\nu $ $=$ $(-1,$ $0)$ on 
$\{$ $(0,y)$ $|$ $t_{1}$ $<$ $y$ $<$ $t_{2}\}.$ It follows from (\ref{4.18})
that (recall that $N^{\perp }$ $=$ $(\sin \theta ,$ $-\cos \theta ))$%
\begin{eqnarray}
&&\int_{y^{t_{1}}(x)}^{y^{t_{2}}(x)}D(x,y)\sin \theta
(x,y)dy+\int_{t_{2}}^{t_{1}}(-D(0,y)\sin \theta (0,y))dy  \label{4.19} \\
&=&\int_{0}^{x}(\int_{y^{t_{1}}(\varsigma )}^{y^{t_{2}}(\varsigma )}rot\vec{F%
}\text{ }d\eta )d\varsigma .  \notag
\end{eqnarray}

\noindent Then by (\ref{4.19}) and (\ref{4.20}) we deduce that%
\begin{equation}
B^{\prime }(x)=\int_{y^{t_{1}}(x)}^{y^{t_{2}}(x)}rot\vec{F}\text{ }%
(x,y)dy=rot\vec{F}(x,\eta ^{\prime })(y^{t_{2}}(x)-y^{t_{1}}(x)).
\label{4.21}
\end{equation}

\noindent and 
\begin{equation}
B(x)=(y^{t_{2}}(x)-y^{t_{1}}(x))D(x,\eta )\sin \theta (x,\eta )  \label{4.22}
\end{equation}

\noindent for $\eta ^{\prime }$ $=$ $\eta ^{\prime }(x,$ $t_{1},$ $t_{2}),$ $%
\eta $ $=$ $\eta (x,$ $t_{1},$ $t_{2})$ such that $y^{t_{1}}(x)$ $<$ $\eta ,$
$\eta ^{\prime }$ $<$ $y^{t_{2}}(x)$ by the mean value theorem. By (\ref%
{4.21}) and (\ref{4.22}) we obtain%
\begin{equation}
\frac{d\log B(x)}{dx}=\frac{B^{\prime }(x)}{B(x)}=\frac{rot\vec{F}(x,\eta
^{\prime })}{D(x,\eta )\sin \theta (x,\eta )}  \label{4.23}
\end{equation}

\noindent (noting that $\sin \theta $ is close to $1$ near $p_{0}$ where $%
\theta $ equals $\frac{\pi }{2}$ by assumption). Integrating both sides of (%
\ref{4.23}) gives%
\begin{equation}
\frac{B(x)}{B(0)}=\exp \int_{0}^{x}\frac{rot\vec{F}(\varsigma ,\eta ^{\prime
})}{D(\varsigma ,\eta )\sin \theta (\varsigma ,\eta )}d\varsigma .
\label{4.24}
\end{equation}

\noindent Observe that $B(0)$ $=$ $(t_{2}$ $-$ $t_{1})$ $D(0,$ $\eta
(0,t_{1},t_{2}))$ $\sin \theta (0,$ $\eta (0,t_{1},t_{2}))$ by (\ref{4.22}).
It then follows from (\ref{4.24}) that%
\begin{eqnarray}
&&\frac{t_{2}-t_{1}}{y^{t_{2}}(x)-y^{t_{1}}(x)}  \label{4.25} \\
&=&\frac{D(x,\eta (x))\sin \theta (x,\eta (x))}{D(0,\eta (0))\sin \theta
(0,\eta (0))}\exp (-\int_{0}^{x}\frac{rot\vec{F}(\varsigma ,\eta ^{\prime
}(\varsigma ))}{D(\varsigma ,\eta (\varsigma ))\sin \theta (\varsigma ,\eta
(\varsigma ))}d\varsigma ).  \notag
\end{eqnarray}

\noindent We have omitted the dependence of $t_{1}$ and $t_{2}$ for the
expression of $\eta $ and $\eta ^{\prime }$ in (\ref{4.25}). Combining (\ref%
{4.3}) with $s$ replaced by $t$ and (\ref{4.25}) and taking the limit $y_{2}$
$\rightarrow $ $y_{1},$ we finally obtain%
\begin{equation}
\frac{\partial t}{\partial y}(x,y_{1})=\frac{D(x,y_{1})\sin \theta (x,y_{1})%
}{D(0,t_{1})\sin \theta (0,t_{1})}\exp (-\int_{0}^{x}\frac{rot\vec{F}%
(\varsigma ,y^{t_{1}}(\varsigma ))}{D(\varsigma ,y^{t_{1}}(\varsigma ))\sin
\theta (\varsigma ,y^{t_{1}}(\varsigma ))}d\varsigma ).  \label{4.26}
\end{equation}

\noindent Here we have used Lebesgue's Dominated Convergence Theorem since $%
\frac{rot\vec{F}(\varsigma ,\eta ^{\prime }(\varsigma ))}{D(\varsigma ,\eta
(\varsigma ))\sin \theta (\varsigma ,\eta (\varsigma ))}$ is uniformly
bounded and both $\eta ^{\prime }$ and $\eta $ converge to $y^{t_{1}}$
pointwise$.$ Let $\sigma $ denote the arc length (unit-speed) parameter of
the characteristic curve $\{t=t_{1}\}.$ By (\ref{4.17}), the definition of $%
t,$ we have $\frac{\partial t}{\partial \sigma }$ $=$ $0.$ Note that $%
N^{\perp }d\sigma $ $=$ $(dx,$ $dy)$ along $\{t$ $=$ $t_{1}\}$ $=$ $\{$ $(x,$
$y^{t_{1}}(x))$ $\}.$ It then follows that%
\begin{eqnarray}
\frac{\partial t}{\partial x}(x,y_{1}) &=&-\frac{\partial t}{\partial y}%
(x,y_{1})\frac{dy^{t_{1}}(x)}{dx}=-\frac{\partial t}{\partial y}(x,y_{1})(-%
\frac{\cos \theta (x,y_{1})}{\sin \theta (x,y_{1})})  \label{4.27} \\
&=&g(x,y_{1})D(x,y_{1})\cos \theta (x,y_{1})  \notag
\end{eqnarray}

\noindent by (\ref{4.26}), where%
\begin{equation}
g(x,y_{1})\equiv \frac{1}{D(0,t_{1})\sin \theta (0,t_{1})}\exp (-\int_{0}^{x}%
\frac{rot\vec{F}(\varsigma ,y^{t_{1}}(\varsigma ))}{D(\varsigma
,y^{t_{1}}(\varsigma ))\sin \theta (\varsigma ,y^{t_{1}}(\varsigma ))}%
d\varsigma ).  \label{4.28}
\end{equation}

\noindent Suppose $\sigma $ $=$ $0$ at $(0,$ $t_{1})$ and $\sigma $ $=$ $l$
at $(x,y_{1})$. Recall that $N^{\perp }d\sigma $ $=$ $(dx,$ $dy)$ and hence $%
\sin \theta (\varsigma ,y^{t_{1}}(\varsigma ))$ $d\sigma $ $=$ $d\varsigma .$
So we can rewrite (\ref{4.28}) as%
\begin{equation}
g(x,y_{1})=g(0,t_{1})\exp (-\int_{0}^{l}\frac{rot\vec{F}(\Gamma
_{(0,t_{1})}(\sigma ))}{D(\Gamma _{(0,t_{1})}(\sigma ))}d\sigma )
\label{4.29}
\end{equation}

\noindent where $\Gamma _{(0,t_{1})}$ denotes the characteristic curve $\{t$ 
$=$ $t_{1}\}$ from $(0,t_{1})$ to $(x,y_{1}),$ parametrized by $\sigma .$
Since $t$ $\in $ $C^{0},$ $\frac{rot\vec{F}(\varsigma ,y^{t_{1}}(\varsigma ))%
}{D(\varsigma ,y^{t_{1}}(\varsigma ))\sin \theta (\varsigma
,y^{t_{1}}(\varsigma ))}$ is uniformly bounded and $y^{t_{2}}$ converges to $%
y^{t_{1}}$ pointwise as $y_{2}$ $\rightarrow $ $y_{1},$ we can apply
Lebesgue's Dominated Convergence Theorem to (\ref{4.28}) and conclude that $%
g $ is continuous in $y_{1}$. On the other hand, $g$ is continuous along the
characteristic curve in view of (\ref{4.29}). Together we have $g$ $\in $ $%
C^{0}$ near $p_{0}$ since the characteristic curves are transverse to the $y$%
-axes $\{x$ $=$ $c,$ constants\} near $p_{0}.$ From (\ref{4.27}), (\ref{4.26}%
), and (\ref{4.28}), we have proved $t$ $\in $ $C^{1}$ and

\begin{equation*}
\nabla t=gD(\cos \theta ,\sin \theta )=gDN
\end{equation*}

\noindent with $g$ $>$ $0$ and $g$ $\in $ $C^{0}$ near $p_{0}.$ Moreover,
recall that $(N^{\perp }g)(x,y_{1})$ $\equiv $ $dg(\Gamma
_{(0,t_{1})}(\sigma ))/d\sigma $ at $\sigma $ $=$ $l$ (if exists)$.$ Now (%
\ref{4.16}) easily follows from (\ref{4.29}).

\endproof%

\bigskip

\proof
\textbf{(of Theorem C) }The\textbf{\ }existence of $s,t$ and (\ref{1.8}), (%
\ref{1.9}) follow from Theorem 4.1 and Theorem 4.2 in view of Theorem 3.2.
By (\ref{1.8}) and (\ref{1.9}) we learn that the Jacobian of $\Psi $ does
not vanish. So by the inverse function theorem $\Psi $ is a $C^{1}$
diffeomorphism near $p_{0}$ and (\ref{1.10}) follows from (\ref{1.8}).

\endproof%

\bigskip

\bigskip 
\proof
\textbf{(of Corollary C.1) }Observe that by (\ref{1.5}) we have\ 
\begin{equation}
y^{\prime }(x)=-\frac{\cos \theta }{\sin \theta }  \label{4.30}
\end{equation}

\noindent (Recall that we have assumed that the characteristic curve $\Gamma 
$ is a piece of a $C^{1}$ smooth graph $(x,$ $y(x))$ so that $\sin \theta $ $%
>$ $0$ without loss of generality$)$. Since $\sin \theta $ $>$ $0,$ we
obtain from (\ref{4.30}) that 
\begin{equation}
\frac{y^{\prime }(x)}{\sqrt{1+(y^{\prime }(x))^{2}}}=-\cos \theta .
\label{4.31}
\end{equation}%
\noindent In the proof of Theorem A, we have shown that the left-hand side
of (\ref{4.31}) is $C^{1}$ smooth in $x.$ It follows that $\theta $ is also $%
C^{1}$ smooth in $x$ and by (\ref{2.7}) we obtain%
\begin{equation}
(\sin \theta )\theta _{x}=-H.  \label{4.32}
\end{equation}

\noindent Recall that $\sigma $ denotes the arc-length parameter. Along $%
\Gamma $ we have%
\begin{equation}
\frac{d\theta }{d\sigma }=\theta _{x}\frac{dx}{d\sigma }=\theta _{x}\sin
\theta  \label{4.33}
\end{equation}%
\bigskip

\noindent in view of (\ref{1.5}). From (\ref{1.10}) we have $(\Psi ^{\ast })$
$ds$ $=$ $fd\sigma $ along $\Gamma $ which is defined by $t$ $=$ constant.
Hence we can compute%
\begin{equation*}
\theta _{s}=\frac{1}{f}\frac{d\theta }{d\sigma }=-\frac{H}{f}
\end{equation*}

\noindent by (\ref{4.33}) and (\ref{4.32}).

\endproof%

\bigskip

\section{Regularity of $\protect\theta $ and characteristic and seed curves}

In Section 4, for a $C^{1}$ smooth solution $u$ to (\ref{1.2}) with $\vec{F}$
$\in $ $C^{1}(\Omega ),$ $H$ $\in $ $C^{0}(\Omega )$, and $\Omega $
nonsingular, we can find locally $C^{1}$ smooth characteristic coordinates $%
s,$ $t,$ that is, a local coordinate change $\Psi :$ $(x,$ $y)\rightarrow $ $%
(s,$ $t)$ which is a $C^{1}$ smooth diffeomorphism such that $\{s$ $=$
constants$\}$ are seed curves and $\{t$ $=$ constants$\}$ are characteristic
curves. In this section we will prove that if $N(H)$ exists and is in $C^{0}$%
, then $\theta $ is $C^{1}$ smooth with respect to $x,$ $y$ coordinates.
This is equivalent to proving that $\theta $ is $C^{1}$ smooth with respect
to $s,$ $t$ coordinates since $\Psi $ is a $C^{1}$ smooth diffeomorphism. So
we consider only $s,$ $t$ coordinates throughout this section.

\bigskip

\textbf{Definition 5.1.} \textit{Let }$\Omega _{1}$\textit{\ and }$\Omega $%
\textit{\ be domains of }$R^{2}$\textit{\ such that }$\Omega _{1}$\textit{\ }%
$\subset \subset $\textit{\ }$\Omega .$\textit{\ We call }$s,$\textit{\ }$t$%
\textit{\ }$C^{1}$\textit{\ coordinates of }$\bar{\Omega}_{1}$\textit{\ if
there exists a domain }$\Omega _{2}$\textit{\ such that }$\Omega _{1}$%
\textit{\ }$\subset \subset $\textit{\ }$\Omega _{2}$\textit{\ }$\subset $%
\textit{\ }$\Omega $\textit{\ and }$s,$\textit{\ }$t$\textit{\ are }$C^{1}$%
\textit{\ coordinates of }$\Omega _{2},$ \textit{i.e., the coordinate change 
}$\Psi :$\textit{\ }$(x,$\textit{\ }$y)\rightarrow $\textit{\ }$(s,$\textit{%
\ }$t)$\textit{\ is a }$C^{1}$\textit{\ smooth diffeomorphism onto }$\Omega
_{2}.$

\bigskip

\textbf{Lemma 5.1.} \textit{Let }$\Omega _{1}$\textit{\ and }$\Omega $%
\textit{\ be domains of }$R^{2}$\textit{\ such that }$\Omega _{1}$\textit{\ }%
$\subset \subset $\textit{\ }$\Omega .$\textit{\ Let }$s,$\textit{\ }$t$%
\textit{\ be }$C^{1}$\textit{\ coordinates of }$\bar{\Omega}_{1}$\textit{\
and }$\Omega _{1}$\textit{\ }$=$\textit{\ }$(0,$\textit{\ }$\tilde{s})$%
\textit{\ }$\times $\textit{\ }$(0,$\textit{\ }$\tilde{t})$\textit{\ for
some }$\tilde{s},$\textit{\ }$\tilde{t}$\textit{\ }$>$\textit{\ }$0.$\textit{%
\ Suppose }$h$\textit{\ }$\in $\textit{\ }$C^{0}(\Omega ),$\textit{\ and for
points in }$\bar{\Omega}_{1},$\textit{\ }$h$\textit{\ is }$C^{1}$\textit{\
smooth in }$s$\textit{\ and }$h_{s}$\textit{\ }$\equiv $\textit{\ }$\frac{%
\partial h}{\partial s}$\textit{\ is }$C^{1}$\textit{\ smooth in }$t.$%
\textit{\ Let }$(s_{2},$\textit{\ }$t_{2}),$\textit{\ }$(s_{1},$\textit{\ }$%
t_{2}),$\textit{\ }$(s_{2},$\textit{\ }$t_{1}),$\textit{\ }$(s_{1},$\textit{%
\ }$t_{1})$\textit{\ }$\in $\textit{\ }$\bar{\Omega}_{1}$\textit{\ with }$%
t_{2}$\textit{\ }$\neq $\textit{\ }$t_{1}.$\textit{\ }

\textit{(a) Let }$M$\textit{\ }$>$\textit{\ }$0$\textit{\ be a given
constant. Let }$s_{0}$\textit{\ }$=$\textit{\ }$\min \{\tilde{s},$\textit{\ }%
$\frac{M}{2}(\max_{\bar{\Omega}_{1}}|(h_{s})_{t}|+1)^{-1}\}.$\textit{\ Then
if }$\frac{h(s_{1},t_{2})-h(s_{1},t_{1})}{t_{2}-t_{1}}$\textit{\ }$>$\textit{%
\ }$M$\textit{\ }$>$\textit{\ }$0$\textit{(or }$<$\textit{\ }$-M,$\textit{\
resp.), then }$\frac{h(s_{2},t_{2})-h(s_{2},t_{1})}{t_{2}-t_{1}}$\textit{\ }$%
>$\textit{\ }$\frac{M}{2}$\textit{\ (}$<$\textit{\ }$-\frac{M}{2},$\textit{\
resp.}$)$\textit{\ for }$|s_{2}$\textit{\ }$-$\textit{\ }$s_{1}|$\textit{\ }$%
\leq $\textit{\ }$s_{0}.$

\textit{(b) Given }$s_{0}^{\prime },$\textit{\ }$0$\textit{\ }$\leq $\textit{%
\ }$s_{0}^{\prime }$\textit{\ }$\leq $\textit{\ }$\tilde{s}.$\textit{\ Let }$%
M_{0}$\textit{\ }$=$\textit{\ }$2s_{0}^{\prime }\max_{\bar{\Omega}%
_{1}}|(h_{s})_{t}|.$\textit{\ Then for any }$M$\textit{\ }$\geq $\textit{\ }$%
M_{0},$\textit{\ }$|s_{2}$\textit{\ }$-$\textit{\ }$s_{1}|$\textit{\ }$\leq $%
\textit{\ }$s_{0}^{\prime }$\textit{, if }$\frac{%
h(s_{1},t_{2})-h(s_{1},t_{1})}{t_{2}-t_{1}}$\textit{\ }$>$\textit{\ }$M$%
\textit{\ (or }$<$\textit{\ }$-M,$\textit{\ resp.), then }$\frac{%
h(s_{2},t_{2})-h(s_{2},t_{1})}{t_{2}-t_{1}}$\textit{\ }$>$\textit{\ }$\frac{M%
}{2}$\textit{\ (}$<$\textit{\ }$-\frac{M}{2},$\textit{\ resp.}$).$

\bigskip

Recall that by $h$\ being $C^{1}$\ smooth in $s$ we mean that $h_{s}$\ $%
\equiv $\ $\frac{\partial h}{\partial s}$ exists and is continuous. Since we
only assume that $h$\ is $C^{1}$\ smooth in $s,$ we need to consider the
behavior of $h$ with respect to $t$ in order to prove $h$ $\in $ $C^{1}$ for
later applications. Therefore we study the properties of difference quotient 
$\frac{h(s,t_{2})-h(s,t_{1})}{t_{2}-t_{1}}.$

\bigskip

\proof
We can write%
\begin{equation}
h(s_{2},t)-h(s_{1},t)=\int_{s_{1}}^{s_{2}}h_{s}(s,t)ds  \label{5.3}
\end{equation}

\noindent by the fundamental theorem of calculus. With $t$ $=$ $t_{2}$ and $%
t_{1}$ in (\ref{5.3}) respectively, we then substract one resulting formula
from the other to get%
\begin{equation}
\frac{h(s_{2},t_{2})-h(s_{2},t_{1})}{t_{2}-t_{1}}-\frac{%
h(s_{1},t_{2})-h(s_{1},t_{1})}{t_{2}-t_{1}}=\int_{s_{1}}^{s_{2}}\frac{%
h_{s}(s,t_{2})-h_{s}(s,t_{1})}{t_{2}-t_{1}}ds.  \label{5.4}
\end{equation}

\noindent By the mean value theorem we can find $t^{\prime }$ $=$ $t^{\prime
}(s)$ between $t_{1}$ and $t_{2}$ such that%
\begin{equation}
\frac{h_{s}(s,t_{2})-h_{s}(s,t_{1})}{t_{2}-t_{1}}=(h_{s})_{t}(s,t^{\prime })
\label{5.5}
\end{equation}

\noindent since $h_{s}$ is $C^{1}$ smooth in $t$ by assumption. Substituting
(\ref{5.5}) into (\ref{5.4}), we obtain%
\begin{equation}
\frac{h(s_{2},t_{2})-h(s_{2},t_{1})}{t_{2}-t_{1}}-\frac{%
h(s_{1},t_{2})-h(s_{1},t_{1})}{t_{2}-t_{1}}%
=\int_{s_{1}}^{s_{2}}(h_{s})_{t}(s,t^{\prime }(s))ds.  \label{5.6}
\end{equation}

\noindent From (\ref{5.6}) we can easily deduce (a) and (b).

\endproof%

\bigskip

\textbf{Lemma 5.2. }\textit{Suppose we have the same situation as in Lemma
5.1. }

\textit{(a) Given }$\varepsilon $\textit{\ }$>$\textit{\ }$0$\textit{\ we
can find }$\delta _{0}$\textit{\ }$=$\textit{\ }$\delta _{0}((h_{s})_{t}|_{%
\bar{\Omega}_{1}},\varepsilon )$\textit{\ }$>$\textit{\ }$0$\textit{\ such
that for }$p,$\textit{\ }$q$\textit{\ }$\in $\textit{\ }$\bar{\Omega}_{1},$%
\textit{\ }$|p-q|$\textit{\ }$<$\textit{\ }$\delta _{0},$\textit{\ there
holds }%
\begin{equation}
|(h_{s})_{t}(p)-(h_{s})_{t}(q)|<\varepsilon .  \label{5.2.0}
\end{equation}

\textit{(b) Let }$(s_{2},$\textit{\ }$t_{3}),$\textit{\ }$(s_{1},$\textit{\ }%
$t_{3}),$\textit{\ }$(s_{2},$\textit{\ }$t_{2}),$\textit{\ }$(s_{1},$\textit{%
\ }$t_{2}),$\textit{\ }$(s_{2},$\textit{\ }$t_{1}),$\textit{\ }$(s_{1},$%
\textit{\ }$t_{1})$\textit{\ }$\in $\textit{\ }$\bar{\Omega}_{1}$\textit{\
with }$t_{3}$\textit{\ }$\neq $\textit{\ }$t_{1},$\textit{\ }$t_{2}$\textit{%
\ }$\neq $\textit{\ }$t_{1}.$\textit{\ Given }$\varepsilon $\textit{\ }$>$%
\textit{\ }$0.$\textit{\ Then for }$|t_{3}-t_{1}|$\textit{\ }$+$\textit{\ }$%
|t_{2}-t_{1}|$\textit{\ }$<$\textit{\ }$\delta _{0}$\textit{\ (as in (a)),
there holds}%
\begin{eqnarray}
&&|[\frac{h(s_{2},t_{3})-h(s_{2},t_{1})}{t_{3}-t_{1}}-\frac{%
h(s_{2},t_{2})-h(s_{2},t_{1})}{t_{2}-t_{1}}]  \label{5.2} \\
&&-[\frac{h(s_{1},t_{3})-h(s_{1},t_{1})}{t_{3}-t_{1}}-\frac{%
h(s_{1},t_{2})-h(s_{1},t_{1})}{t_{2}-t_{1}}]|  \notag \\
&\leq &|s_{2}-s_{1}|\varepsilon .  \notag
\end{eqnarray}

\bigskip

\proof
(a) follows by observing that $(h_{s})_{t}$ is uniformly continuous on $\bar{%
\Omega}_{1}$ since it is continuous on $\bar{\Omega}_{1}$ and $\bar{\Omega}%
_{1}$ is compact. For the proof of (b), following the proof of Lemma 5.1,
similarly to (\ref{5.6}), we can find $t^{\prime \prime }$ between $t_{1}$
and $t_{3}$ such that 
\begin{equation}
\frac{h(s_{2},t_{3})-h(s_{2},t_{1})}{t_{3}-t_{1}}-\frac{%
h(s_{1},t_{3})-h(s_{1},t_{1})}{t_{3}-t_{1}}%
=\int_{s_{1}}^{s_{2}}(h_{s})_{t}(s,t^{\prime \prime })ds.  \label{5.7}
\end{equation}

\noindent Substracting (\ref{5.6}) from (\ref{5.7}), we can then estimate%
\begin{eqnarray*}
&&|[\frac{h(s_{2},t_{3})-h(s_{2},t_{1})}{t_{3}-t_{1}}-\frac{%
h(s_{1},t_{3})-h(s_{1},t_{1})}{t_{3}-t_{1}}] \\
&&-[\frac{h(s_{2},t_{2})-h(s_{2},t_{1})}{t_{2}-t_{1}}-\frac{%
h(s_{1},t_{2})-h(s_{1},t_{1})}{t_{2}-t_{1}}]| \\
&\leq &\int_{s_{1}}^{s_{2}}|(h_{s})_{t}(s,t^{\prime \prime
})-(h_{s})_{t}(s,t^{\prime })|ds \\
&\leq &|s_{2}-s_{1}|\varepsilon
\end{eqnarray*}

\noindent by (\ref{5.2.0}) since $|(s,t^{\prime \prime })$ $-$ $(s,t^{\prime
})|$ $=$ $|t^{\prime \prime }$ $-$ $t^{\prime }|$ $\leq $ $|t_{3}-t_{1}|$%
\textit{\ }$+$\textit{\ }$|t_{2}-t_{1}|$\textit{\ }$\leq $\textit{\ }$\delta
_{0}$.

\endproof%

\bigskip

\textbf{Lemma 5.3.} \textit{Let }$u$\textit{\ }$:$\textit{\ }$\Omega $%
\textit{\ }$\subset $\textit{\ }$R^{2}$\textit{\ }$\rightarrow $\textit{\ }$%
R $\textit{\ be a }$C^{1}$\textit{\ smooth solution to (\ref{1.2}) with }$%
\vec{F}$\textit{\ }$\in $\textit{\ }$C^{1}(\Omega )$\textit{\ and }$H$%
\textit{\ }$\in $\textit{\ }$C^{0}(\Omega )$\textit{\ such that }$\Omega $%
\textit{\ is nonsingular}$.$\textit{\ Let }$p_{0}$\textit{\ }$\in $\textit{\ 
}$\Omega $ \textit{and }$\theta $\textit{\ (defined locally) }$\in $\textit{%
\ }$C^{0}$\textit{\ such that }$(\cos \theta ,$\textit{\ }$\sin \theta )$%
\textit{\ }$=$\textit{\ }$N_{\vec{F}}^{u}$\textit{\ with }$\theta (p_{0})$%
\textit{\ }$=$\textit{\ }$\frac{\pi }{2}.$\textit{\ Then }

\textit{(a) there exists a domain }$\Omega _{1}$\textit{\ such that }$p_{0}$%
\textit{\ }$\in $\textit{\ }$\Omega _{1}$\textit{\ }$\subset \subset $%
\textit{\ }$\Omega $\textit{, }$\bar{\Omega}_{1}$\textit{\ has }$C^{1}$%
\textit{\ coordinates }$s,$\textit{\ }$t,$\textit{\ }$\Omega _{1}$\textit{\ }%
$=$\textit{\ }$(0,\tilde{s})$\textit{\ }$\times $\textit{\ }$(0,\tilde{t})$%
\textit{\ in }$s,$\textit{\ }$t$\textit{\ coordinates, and }$|\theta $%
\textit{\ }$-$\textit{\ }$\frac{\pi }{2}|$\textit{\ }$<<$\textit{\ }$1$%
\textit{\ in }$\bar{\Omega}_{1}.$

\textit{(b) Take }$s,$\textit{\ }$t$\textit{\ coordinates as in (a). Suppose 
}$H$\textit{\ is }$C^{1}$\textit{\ smooth in }$t.$\textit{\ Then there
exists a constant }$M$\textit{\ }$=$\textit{\ }$M(\tilde{s},$\textit{\ }$%
\max_{\bar{\Omega}_{1}}|y_{t}|,$\textit{\ }$\max_{\bar{\Omega}_{1}}|f|,$%
\textit{\ }$\max_{\bar{\Omega}_{1}}|\frac{f_{t}}{f^{2}}|,$\textit{\ }$\max_{%
\bar{\Omega}_{1}}|(\theta _{s})_{t}|)$\textit{\ }$>$\textit{\ }$0$\textit{\
such that }%
\begin{equation}
\mid \frac{\theta (s,t_{2})-\theta (s,t_{1})}{t_{2}-t_{1}}\mid \leq M
\label{5.7.1}
\end{equation}%
\textit{\noindent for any }$(s,$\textit{\ }$t_{2}),$\textit{\ }$(s,$\textit{%
\ }$t_{1})$\textit{\ }$\in $\textit{\ }$\bar{\Omega}_{1}$\textit{\ and }$%
t_{2}$\textit{\ }$\neq $\textit{\ }$t_{1}.$

\textit{\bigskip }

\proof%
\textbf{\ }(a) follows\textbf{\ }from Theorem C and $\theta $ $\in $ $C^{0}.$
To prove (b), by (\ref{1.10}) we first rewrite the equations (\ref{1.5}) of
the characteristic curves in $(s,$ $t)$ coordinates:

\begin{eqnarray}
\frac{dx(s,t)}{ds} &=&\frac{\sin \theta (s,t)}{f(s,t)},  \label{5.8} \\
\frac{dy(s,t)}{ds} &=&-\frac{\cos \theta (s,t)}{f(s,t)}.  \notag
\end{eqnarray}

\noindent From the second formula of (\ref{5.8}) we have, for $0$ $\leq $ $%
s_{1},$ $s_{2}$ $\leq $ $\tilde{s},$ $0$ $\leq $ $t_{1},$ $t_{2}$ $\leq $ $%
\tilde{t}$ and $t_{2}$\textit{\ }$\neq $\textit{\ }$t_{1},$%
\begin{eqnarray}
&&\frac{y(s_{2},t_{2})-y(s_{2},t_{1})}{t_{2}-t_{1}}-\frac{%
y(s_{1},t_{2})-y(s_{1},t_{1})}{t_{2}-t_{1}}  \label{5.10} \\
&=&\frac{y(s_{2},t_{2})-y(s_{1},t_{2})}{t_{2}-t_{1}}-\frac{%
y(s_{2},t_{1})-y(s_{1},t_{1})}{t_{2}-t_{1}}  \notag \\
&=&\int_{s_{1}}^{s_{2}}\frac{1}{t_{2}-t_{1}}[\frac{\cos \theta (s,t_{1})}{%
f(s,t_{1})}-\frac{\cos \theta (s,t_{2})}{f(s,t_{2})}]ds  \notag \\
&=&\int_{s_{1}}^{s_{2}}(A+B)ds.  \notag
\end{eqnarray}

\noindent where 
\begin{eqnarray}
A &=&\frac{-1}{f(s,t_{2})}\frac{\cos \theta (s,t_{2})-\cos \theta (s,t_{1})}{%
t_{2}-t_{1}},  \label{5.11} \\
B &=&-\frac{\cos \theta (s,t_{1})}{t_{2}-t_{1}}[\frac{1}{f(s,t_{2})}-\frac{1%
}{f(s,t_{1})}].  \notag
\end{eqnarray}

Observe that 
\begin{equation}
A=\frac{\sin \theta ^{\prime }}{f(s,t_{2})}\frac{\theta (s,t_{2})-\theta
(s,t_{1})}{t_{2}-t_{1}}  \label{5.12}
\end{equation}

\noindent where $\theta ^{\prime }$ is a number between $\theta (s,t_{2})$
and $\theta (s,t_{1})$ and

\begin{equation}
B=\frac{\cos \theta (s,t_{1})}{f^{2}(s,t^{\prime })}\frac{\partial f}{%
\partial t}(s,t^{\prime })  \label{5.13}
\end{equation}

\noindent where $t^{\prime }$ is a number between $t_{2}$ and $t_{1}$ by the
mean value theorem (noting that $f$ $>$ $0$ is $C^{1}$ smooth in $t$ by
Theorem C). Since $|\theta $\textit{\ }$-$\textit{\ }$\frac{\pi }{2}|$%
\textit{\ }$<<$\textit{\ }$1$\textit{\ }in\textit{\ }$\bar{\Omega}_{1},$ we
have%
\begin{equation}
\frac{\sin \theta ^{\prime }}{f(s,t_{2})}\geq C_{1}\equiv \frac{1}{2\max_{%
\bar{\Omega}_{1}}|f|},\text{ }\mid B\mid \leq C_{2}\equiv \max_{\bar{\Omega}%
_{1}}|\frac{f_{t}}{f^{2}}|  \label{5.14}
\end{equation}

\noindent in $\bar{\Omega}_{1}.$

On the other hand, we have%
\begin{eqnarray}
&&\left\vert \frac{y(s_{2},t_{2})-y(s_{2},t_{1})}{t_{2}-t_{1}}-\frac{%
y(s_{1},t_{2})-y(s_{1},t_{1})}{t_{2}-t_{1}}\right\vert  \label{5.15} \\
&\leq &2\max_{\bar{\Omega}_{1}}|y_{t}|\leq C_{3}\equiv 2\max_{\bar{\Omega}%
_{1}}|y_{t}|+1  \notag
\end{eqnarray}

\noindent since $y$ is $C^{1}$ smooth on $\bar{\Omega}_{1}.$ Since $\theta
_{s}$ $=$ $-\frac{H}{f}$ (see Corollary C.1)$,$ $(\theta _{s})_{t}$ $=$ $-(%
\frac{H}{f})_{t}$ $\in $ $C^{0}$ by the assumption on $H$ and $f$ being $%
C^{1}$ in $t$ due to Theorem C. So the assumptions of Lemma 5.1 are
satisfied for $h$ $=$ $\theta .$ Take $s_{0}^{\prime }$ $=$ $\tilde{s}$ in
Lemma 5.1 (b). Since $0$ $\leq $ $s_{1},$ $s_{2}$ $\leq $ $\tilde{s},$ we
have $|s_{2}$ $-$ $s_{1}|$ $\leq $ $\tilde{s}.$ Let $M$ $=$ $\max \{M_{0}$ $%
\equiv $ $2\tilde{s}\max_{\bar{\Omega}_{1}}|(\theta _{s})_{t}|,$ $\frac{%
2(2C_{3}+\tilde{s}C_{2})}{\tilde{s}C_{1}}\}$. Suppose (\ref{5.7.1}) fails to
hold. Then we can find $(s_{1},$ $t_{2}),$ $(s_{1},$ $t_{1})$ $\in $ $\bar{%
\Omega}_{1},$ $t_{2}$\textit{\ }$\neq $\textit{\ }$t_{1}$, such that 
\begin{equation*}
\frac{\theta (s_{1},t_{2})-\theta (s_{1},t_{1})}{t_{2}-t_{1}}>M\text{ }(%
\text{or }<-M,\text{ resp.}).
\end{equation*}

\noindent Applying Lemma 5.1(b) to $h$ $=$ $\theta $ gives%
\begin{equation}
\frac{\theta (s_{2},t_{2})-\theta (s_{2},t_{1})}{t_{2}-t_{1}}>\frac{M}{2}%
\text{ }(\text{or }<-\frac{M}{2},\text{ resp.}).  \label{5.15.0}
\end{equation}%
\noindent for all $s_{2}$ $\in $ $[0,$\textit{\ }$\tilde{s}].$ From (\ref%
{5.10}), (\ref{5.15}), (\ref{5.14}), and (\ref{5.15.0}), we estimate%
\begin{eqnarray}
C_{3} &\geq &|\frac{y(s_{2}^{\prime },t_{2})-y(s_{2}^{\prime },t_{1})}{%
t_{2}-t_{1}}-\frac{y(s_{1},t_{2})-y(s_{1},t_{1})}{t_{2}-t_{1}}|  \label{5.16}
\\
&=&\left\vert \int_{s_{1}}^{s_{2}^{\prime }}(A+B)ds\right\vert  \notag \\
&>&|s_{2}^{\prime }-s_{1}|(\frac{M}{2}C_{1}-C_{2})\geq \frac{\tilde{s}}{2}(%
\frac{M}{2}C_{1}-C_{2})\geq C_{3}  \notag
\end{eqnarray}

\noindent for some $s_{2}^{\prime },$ $0$ $\leq $ $s_{2}^{\prime }$ $\leq $ $%
\tilde{s},$ satisfying $|s_{2}^{\prime }$ $-$ $s_{1}|$ $\geq $ $\frac{\tilde{%
s}}{2}.$ We have reached a contradiction. Therefore (\ref{5.7.1}) holds.

\endproof%

\bigskip

The following lemma should be a basic fact in calculus. For completeness we
give a proof.

\bigskip

\textbf{Lemma 5.4.} \textit{Let }$w$\textit{\ be a real }$C^{1}$\textit{\
smooth function of }$(s,$\textit{\ }$t)$\textit{\ }$\in $\textit{\ }$\Omega $%
\textit{\ }$\subset $\textit{\ }$R^{2}.$\textit{\ Suppose further that (}$%
w_{s})_{t}$\textit{\ exists and is continuous. Then }$(w_{t})_{s}$\textit{\
exists and equals (}$w_{s})_{t}$\textit{\ (hence is continuous)}.

\bigskip

\proof
Take $(s_{0},$ $t_{0})$ $\in $ $\Omega .$ Then there exists $r_{0}$ $>$ $0$
such that $B_{r_{0}}((s_{0},$ $t_{0}))$ $\subset $ $\Omega .$ For $(s_{1},$ $%
t_{1})$ $\in $ $B_{r_{0}}((s_{0},$ $t_{0}))$ we write%
\begin{eqnarray}
w(s_{1},t_{1})-w(s_{0},t_{1}) &=&\int_{s_{0}}^{s_{1}}w_{s}(s,t_{1})ds
\label{5.17.1} \\
w(s_{1},t_{0})-w(s_{0},t_{0}) &=&\int_{s_{0}}^{s_{1}}w_{s}(s,t_{0})ds. 
\notag
\end{eqnarray}

\noindent Take the difference of the two equalities in (\ref{5.17.1})
divided by $t_{1}$ $-$ $t_{0}.$ Let $t_{1}$ approach $t_{0}$ in the
resulting formula to get%
\begin{equation}
w_{t}(s_{1},t_{0})-w_{t}(s_{0},t_{0})=%
\int_{s_{0}}^{s_{1}}(w_{s})_{t}(s,t_{0})ds.  \label{5.17.2}
\end{equation}

\noindent Now dividing (\ref{5.17.2}) by $s_{1}$ $-$ $s_{0}$ and applying
the mean value theorem to $(w_{s})_{t}$, we obtain%
\begin{eqnarray}
\frac{w_{t}(s_{1},t_{0})-w_{t}(s_{0},t_{0})}{s_{1}-s_{0}} &=&\frac{1}{%
s_{1}-s_{0}}\int_{s_{0}}^{s_{1}}(w_{s})_{t}(s,t_{0})ds  \label{5.17.3} \\
&=&(w_{s})_{t}(\tilde{s},t_{0})  \notag
\end{eqnarray}

\noindent for $s_{0}$ $<$ $\tilde{s}$ $<$ $s_{1}.$ Letting $s_{1}$ tend to $%
s_{0}$ in (\ref{5.17.3}), we obtain the existence of $%
(w_{t})_{s}(s_{0},t_{0})$ and 
\begin{equation*}
(w_{t})_{s}(s_{0},t_{0})=(w_{s})_{t}(s_{0},t_{0})
\end{equation*}

\noindent since $(w_{s})_{t}$ is continuous.

\endproof%

\bigskip

\proof
\textbf{(of Theorem D) }Given\textbf{\ }a point $p_{0}$ $\in $ $\Omega ,$ we
may assume $\theta (p_{0})$ $=$ $\frac{\pi }{2}$ without loss of generality.
By Theorem C and $\theta $ (defined locally) $\in $ $C^{0},$ we can choose a
domain $\Omega _{1},$ $p_{0}$ $\in $ $\Omega _{1}$ $\subset \subset $ $%
\Omega ,$ such that $\bar{\Omega}_{1}$\textit{\ }has $C^{1}$\ coordinates $%
s, $\ $t,$\ $\Omega _{1}$\ $=$\ $(0,\tilde{s})$\ $\times $\ $(0,\tilde{t})$\
in $s,$\ $t$\ coordinates, and $|\theta $\ $-$\ $\frac{\pi }{2}|$\ $<<$\ $1$%
\ in $\bar{\Omega}_{1}.$ Since $N$ $=$ $gD\frac{\partial }{\partial t}$ in
view of (\ref{1.10}), we have $H_{t}$ $=$ $\frac{1}{gD}N(H)$ $\in $ $C^{0}$.
It follows from Lemma 5.3 (b) that 
\begin{equation*}
\mid \frac{\theta (s,t_{2})-\theta (s,t_{1})}{t_{2}-t_{1}}\mid \leq M=M(%
\tilde{s},\mathit{\ }\max_{\bar{\Omega}_{1}}|y_{t}|,\mathit{\ }\max_{\bar{%
\Omega}_{1}}|f|,\mathit{\ }\max_{\bar{\Omega}_{1}}|\frac{f_{t}}{f^{2}}|,%
\mathit{\ }\max_{\bar{\Omega}_{1}}|(\theta _{s})_{t}|)
\end{equation*}%
\noindent for any $(s,$\ $t_{2}),$\ $(s,$\ $t_{1})$\ $\in $\ $\bar{\Omega}%
_{1}$\ and $t_{2}$\ $\neq $\ $t_{1}$ (see (\ref{5.7.1}))$.$ We claim that
for $(s_{1},$\ $t),$\ $(s_{1},$\ $t_{1})$\ $\in $\ $\Omega _{1}$\ and $t$\ $%
\neq $\ $t_{1},$%
\begin{equation*}
\frac{\theta (s_{1},t)-\theta (s_{1},t_{1})}{t-t_{1}}
\end{equation*}%
\noindent is a Cauchy sequence as $t$ $\rightarrow $ $t_{1}$ (if so, $\frac{%
\partial \theta }{\partial t}$ exists at $(s_{1},\ t_{1})).$ If not, there
exists $\kappa $ $>$ $0$ such that for any positive integer $n,$ there exist 
$t_{n}^{\prime \prime },$ $t_{n}^{\prime }$ $\in $ $(0,$ $\tilde{t}),$ $%
t_{n}^{\prime \prime }$ $\neq $ $t_{1},$ $t_{n}^{\prime }$ $\neq $ $t_{1}$
satisfying $|t_{n}^{\prime \prime }$ $-$ $t_{1}|$ $+$ $|t_{n}^{\prime }$ $-$ 
$t_{1}|$ $\leq $ $\min (\tilde{t},$ $\frac{1}{n})$ and 
\begin{equation}
\frac{\theta (s_{1},t_{n}^{\prime \prime })-\theta (s_{1},t_{1})}{%
t_{n}^{\prime \prime }-t_{1}}-\frac{\theta (s_{1},t_{n}^{\prime })-\theta
(s_{1},t_{1})}{t_{n}^{\prime }-t_{1}}>\kappa \text{ (or }<-\kappa ,\text{
resp.}).  \label{5.18}
\end{equation}

\noindent Take $h$ $=$ $\theta $ and $\varepsilon $ $=$ $\frac{\kappa }{2%
\tilde{s}}$ in Lemma 5.2. Then there exists $0$ $<$ $\delta _{0}$ $=$ $%
\delta _{0}((\theta _{s})_{t}|_{\bar{\Omega}_{1}},$ $\varepsilon )$ $\leq $ $%
\tilde{t}$ such that 
\begin{eqnarray*}
&&|(\frac{\theta (s_{2},t_{n}^{\prime \prime })-\theta (s_{2},t_{1})}{%
t_{n}^{\prime \prime }-t_{1}}-\frac{\theta (s_{2},t_{n}^{\prime })-\theta
(s_{2},t_{1})}{t_{n}^{\prime }-t_{1}}) \\
&&-(\frac{\theta (s_{1},t_{n}^{\prime \prime })-\theta (s_{1},t_{1})}{%
t_{n}^{\prime \prime }-t_{1}}-\frac{\theta (s_{1},t_{n}^{\prime })-\theta
(s_{1},t_{1})}{t_{n}^{\prime }-t_{1}})| \\
&\leq &|s_{2}-s_{1}|\varepsilon \leq \frac{\kappa }{2}
\end{eqnarray*}%
\noindent for $0$ $\leq $ $s_{2}$ $\leq $ $\tilde{s},$ $|t_{n}^{\prime
\prime }$ $-$ $t_{1}|$ $+$ $|t_{n}^{\prime }$ $-$ $t_{1}|$ $\leq $ $\delta
_{0}$ (see (\ref{5.2})). So if we require $t_{n}^{\prime \prime }$ $\neq $ $%
t_{1},$ $t_{n}^{\prime }$ $\neq $ $t_{1}$ satisfying $|t_{n}^{\prime \prime
} $ $-$ $t_{1}|$ $+$ $|t_{n}^{\prime }$ $-$ $t_{1}|$ $\leq $ $\min (\tilde{t}%
,$ $\frac{1}{n},$ $\delta _{0}),$ we can then estimate%
\begin{eqnarray}
&&\frac{\theta (s_{2},t_{n}^{\prime \prime })-\theta (s_{2},t_{1})}{%
t_{n}^{\prime \prime }-t_{1}}-\frac{\theta (s_{2},t_{n}^{\prime })-\theta
(s_{2},t_{1})}{t_{n}^{\prime }-t_{1}}  \label{5.19} \\
&\geq &\frac{\theta (s_{1},t_{n}^{\prime \prime })-\theta (s_{1},t_{1})}{%
t_{n}^{\prime \prime }-t_{1}}-\frac{\theta (s_{1},t_{n}^{\prime })-\theta
(s_{1},t_{1})}{t_{n}^{\prime }-t_{1}}-\frac{\kappa }{2}  \notag \\
&>&\kappa -\frac{\kappa }{2}=\frac{\kappa }{2}\text{ (}<-\frac{\kappa }{2},%
\text{ resp.})  \notag
\end{eqnarray}

\noindent by (\ref{5.18}) for $0$ $\leq $ $s_{2}$ $\leq $ $\tilde{s}.$ By (%
\ref{5.10}) we have%
\begin{eqnarray}
&&\frac{y(s_{2},t_{n}^{\prime })-y(s_{2},t_{1})}{t_{n}^{\prime }-t_{1}}-%
\frac{y(s_{1},t_{n}^{\prime })-y(s_{1},t_{1})}{t_{n}^{\prime }-t_{1}}
\label{5.19.1} \\
&=&\int_{s_{1}}^{s_{2}}(A_{n}^{\prime }+B_{n}^{\prime })ds  \notag
\end{eqnarray}%
\noindent where 
\begin{eqnarray}
A_{n}^{\prime } &=&\frac{-1}{f(s,t_{n}^{\prime })}\frac{\cos \theta
(s,t_{n}^{\prime })-\cos \theta (s,t_{1})}{t_{n}^{\prime }-t_{1}},
\label{5.19.2} \\
B_{n}^{\prime } &=&-\frac{\cos \theta (s,t_{1})}{t_{n}^{\prime }-t_{1}}[%
\frac{1}{f(s,t_{n}^{\prime })}-\frac{1}{f(s,t_{1})}].  \notag
\end{eqnarray}%
\noindent Since $|B_{n}^{\prime }|$ $\leq $ $\max_{\bar{\Omega}_{1}}|\frac{%
f_{t}}{f^{2}}|$ and $B_{n}^{\prime }$ converges to $\cos \theta (s,t_{1})$ $%
\frac{f_{t}}{f^{2}}(s,t_{1})$ pointwise, we conclude from Lebesgue's
Dominated Convergence Theorem that $\lim_{n\rightarrow \infty
}\int_{s_{1}}^{s_{2}}A_{n}^{\prime }ds$ exists and%
\begin{equation}
\lim_{n\rightarrow \infty }\int_{s_{1}}^{s_{2}}A_{n}^{\prime
}ds=y_{t}(s_{2},t_{1})-y_{t}(s_{1},t_{1})-\int_{s_{1}}^{s_{2}}\cos \theta
(s,t_{1})\frac{f_{t}}{f^{2}}(s,t_{1})ds.  \label{5.20.1}
\end{equation}%
\noindent Similarly with $t_{n}^{\prime }$ replaced by $t_{n}^{\prime \prime
}$ in (\ref{5.19.1}) and the same reasoning, we also obtain%
\begin{equation}
\lim_{n\rightarrow \infty }\int_{s_{1}}^{s_{2}}A_{n}^{\prime \prime
}ds=y_{t}(s_{2},t_{1})-y_{t}(s_{1},t_{1})-\int_{s_{1}}^{s_{2}}\cos \theta
(s,t_{1})\frac{f_{t}}{f^{2}}(s,t_{1})ds  \label{5.20.2}
\end{equation}%
\noindent where 
\begin{equation*}
A_{n}^{\prime \prime }=\frac{-1}{f(s,t_{n}^{\prime \prime })}\frac{\cos
\theta (s,t_{n}^{\prime \prime })-\cos \theta (s,t_{1})}{t_{n}^{\prime
\prime }-t_{1}}.
\end{equation*}%
\noindent Therefore by (\ref{5.20.1}) and (\ref{5.20.2}) we have%
\begin{equation}
\lim_{n\rightarrow \infty }\int_{s_{1}}^{s_{2}}(A_{n}^{\prime \prime
}-A_{n}^{\prime })ds=0.  \label{5.21}
\end{equation}%
\noindent On the other hand, using the mean value theorem, we can find $%
\tilde{\theta}^{\prime }$ ($\tilde{\theta}^{\prime \prime },$ resp.) between 
$\theta (s,t_{n}^{\prime })$ ($\theta (s,t_{n}^{\prime \prime }),$ resp.)
and $\theta (s,t_{1})$ such that 
\begin{eqnarray*}
A_{n}^{\prime } &=&\frac{\sin \tilde{\theta}^{\prime }}{f(s,t_{n}^{\prime })}%
\frac{\theta (s,t_{n}^{\prime })-\theta (s,t_{1})}{t_{n}^{\prime }-t_{1}}, \\
A_{n}^{\prime \prime } &=&\frac{\sin \tilde{\theta}^{\prime \prime }}{%
f(s,t_{n}^{\prime \prime })}\frac{\theta (s,t_{n}^{\prime \prime })-\theta
(s,t_{1})}{t_{n}^{\prime \prime }-t_{1}}.
\end{eqnarray*}

\noindent Hence we can write%
\begin{eqnarray}
A_{n}^{\prime \prime }-A_{n}^{\prime } &=&\frac{\sin \tilde{\theta}^{\prime
\prime }}{f(s,t_{n}^{\prime \prime })}[\frac{\theta (s,t_{n}^{\prime \prime
})-\theta (s,t_{1})}{t_{n}^{\prime \prime }-t_{1}}-\frac{\theta
(s,t_{n}^{\prime })-\theta (s,t_{1})}{t_{n}^{\prime }-t_{1}}]  \label{5.22}
\\
&&+(\frac{\sin \tilde{\theta}^{\prime \prime }}{f(s,t_{n}^{\prime \prime })}-%
\frac{\sin \tilde{\theta}^{\prime }}{f(s,t_{n}^{\prime })})\frac{\theta
(s,t_{n}^{\prime })-\theta (s,t_{1})}{t_{n}^{\prime }-t_{1}}  \notag
\end{eqnarray}

\noindent The second term in the right-hand side of (\ref{5.22}) is
uniformly bounded by Lemma 5.3 (b) and converges to zero pointwise.
Therefore the integral (from $s_{1}$ to $s_{2})$ of this term goes to zero
as $n\rightarrow \infty $ (while $t_{n}^{\prime \prime },$ $t_{n}^{\prime
}\rightarrow $ $t_{1}$ and $\tilde{\theta}^{\prime \prime },$ $\tilde{\theta}%
^{\prime }\rightarrow $ $\theta (s,t_{1}))$ by Lebesgue's Dominated
Convergence Theorem. But the first term in the right-hand side of (\ref{5.22}%
) $\geq $ $\frac{1}{2}(\max_{\bar{\Omega}_{1}}|f|)^{-1}\frac{\kappa }{2}$
(or $\leq $ $-\frac{1}{2}(\max_{\bar{\Omega}_{1}}|f|)^{-1}\frac{\kappa }{2},$
resp.) by (\ref{5.19}). Altogether from (\ref{5.22}) we have the following
estimate%
\begin{equation*}
\mid \lim_{n\rightarrow \infty }\int_{s_{1}}^{s_{2}}(A_{n}^{\prime \prime
}-A_{n}^{\prime })ds\mid \geq |s_{2}-s_{1}|\frac{\kappa }{4}(\max_{\bar{%
\Omega}_{1}}|f|)^{-1}
\end{equation*}%
\noindent which contradicts (\ref{5.21}). Therefore the claim holds, which
implies the existence of $\theta _{t}$ $\equiv $ $\frac{\partial \theta }{%
\partial t}$ at any point $p_{0}.$

For continuity of $\theta _{t},$ we start with (\ref{5.4}) taking $h$ $=$ $%
\theta .$ For $0$ $\leq $ $s_{1},$ $s^{\prime }$ $\leq $ $\tilde{s},$ $0$ $%
\leq $ $t_{1},$ $t^{\prime }$ $\leq $ $\tilde{t},$ $t^{\prime }$ $\neq $ $%
t_{1},$ we have%
\begin{equation}
\frac{\theta (s^{\prime },t^{\prime })-\theta (s^{\prime },t_{1})}{t^{\prime
}-t_{1}}-\frac{\theta (s_{1},t^{\prime })-\theta (s_{1},t_{1})}{t^{\prime
}-t_{1}}=\int_{s_{1}}^{s^{\prime }}\frac{\theta _{s}(s,t^{\prime })-\theta
_{s}(s,t_{1})}{t^{\prime }-t_{1}}ds  \label{5.26}
\end{equation}

\noindent where $\theta _{s}$ $=$ $-\frac{H}{f}$ $\in $ $C^{1}$ in $t$ by (%
\ref{1.10.1}), Theorem C, and the assumption (note that $NH$ $=$ $gD\frac{%
\partial }{\partial t}H)$. Taking the limit $t^{\prime }$ $\rightarrow $ $%
t_{1}$ in (\ref{5.26}), we obtain 
\begin{equation}
\theta _{t}(s^{\prime },t_{1})-\theta
_{t}(s_{1},t_{1})=\int_{s_{1}}^{s^{\prime }}(\theta _{s})_{t}(s,t_{1})ds.
\label{5.27}
\end{equation}

\noindent It follows from (\ref{5.27}) that $\theta _{t}$ is continuous in
the $s$ direction at $(s_{1},t_{1}),$ for all $0$ $\leq $ $s_{1}$ $\leq $ $%
\tilde{s},$ $0$ $\leq $ $t_{1}$ $\leq $ $\tilde{t}.$ So to prove $\theta
_{t} $ is continuous at $(s_{1},$ $t_{1}),$ it suffices to show that $\theta
_{t}$ is continuous in the $t$ direction at $(s_{1},$ $t_{1}),$ for all $0$ $%
\leq $ $s_{1}$ $\leq $ $\tilde{s},$ $0$ $\leq $ $t_{1}$ $\leq $ $\tilde{t}.$
Suppose this fails to hold. Then there exists $M^{\prime }$ $>$ $0$ such
that for any positive integer $m,$ we can find $t_{m}^{\prime }$ $\neq $ $%
t_{1}$ satisfying $|t_{m}^{\prime }$ $-$ $t_{1}|$ $<$ $\min \{\tilde{t},$ $%
\frac{1}{m}\}$ and 
\begin{equation}
\theta _{t}(s_{1},t_{m}^{\prime })-\theta _{t}(s_{1},t_{1})>M^{\prime }\text{
(or }<-M^{\prime },\text{ resp.).}  \label{5.27.1}
\end{equation}%
\noindent We take the difference of the formula (\ref{5.27}) with $t_{1}$
replaced by $t_{m}^{\prime }$ and the formula (\ref{5.27}) itself to obtain%
\begin{eqnarray}
\theta _{t}(s^{\prime },t_{m}^{\prime })-\theta _{t}(s^{\prime },t_{1})
&=&\theta _{t}(s_{1},t_{m}^{\prime })-\theta _{t}(s_{1},t_{1})
\label{5.27.2} \\
&&+\int_{s_{1}}^{s^{\prime }}((\theta _{s})_{t}(s,t_{m}^{\prime })-(\theta
_{s})_{t}(s,t_{1}))ds.  \notag
\end{eqnarray}%
\noindent For $|s^{\prime }-s_{1}|$ $\leq $ $\min \{\frac{M^{\prime }}{4}%
(\max_{\bar{\Omega}_{1}}|(\theta _{s})_{t}|)^{-1},$ $\tilde{s}\}$ we then
have%
\begin{equation}
\theta _{t}(s^{\prime },t_{m}^{\prime })-\theta _{t}(s^{\prime },t_{1})>%
\frac{M^{\prime }}{2}\text{ (or }<-\frac{M^{\prime }}{2},\text{ resp.})
\label{5.27.3}
\end{equation}%
\noindent by (\ref{5.27.1}) and (\ref{5.27.2}). On the other hand, we have a
similar formula as (\ref{5.20.1}):%
\begin{equation*}
\lim_{m\rightarrow \infty }\int_{s_{1}}^{s^{\prime }}A_{m}^{\prime
}ds=y_{t}(s^{\prime },t_{1})-y_{t}(s_{1},t_{1})-\int_{s_{1}}^{s^{\prime
}}\cos \theta (s,t_{1})\frac{f_{t}}{f^{2}}(s,t_{1})ds.
\end{equation*}%
\noindent where%
\begin{equation*}
A_{m}^{\prime }=\frac{-1}{f(s,t_{m}^{\prime })}\frac{\cos \theta
(s,t_{m}^{\prime })-\cos \theta (s,t_{1})}{t_{m}^{\prime }-t_{1}}.
\end{equation*}%
\noindent Since $A_{m}^{\prime }$ is uniformly bounded ($\theta $ is
Lipschitzian) and converges pointwise as $m\rightarrow \infty $ ($\theta
_{t} $ exists)$,$ we obtain 
\begin{eqnarray}
y_{t}(s^{\prime },t_{1})-y_{t}(s_{1},t_{1}) &=&\int_{s_{1}}^{s^{\prime }}%
\frac{\sin \theta (s,t_{1})}{f(s,t_{1})}\theta _{t}(s,t_{1})  \label{5.28} \\
&&+\int_{s_{1}}^{s^{\prime }}\cos \theta (s,t_{1})\frac{f_{t}}{f^{2}}%
(s,t_{1})ds  \notag
\end{eqnarray}

\noindent by Lebesgue's dominated convergence theorem. Replacing $t_{1}$ by $%
t_{m}^{\prime }$ in (\ref{5.28}), taking the difference of the resulting
formula and (\ref{5.28}), and letting $m$ $\rightarrow $ $\infty ,$ we
finally obtain%
\begin{equation}
\lim_{m\rightarrow \infty }\int_{s_{1}}^{s^{\prime }}\frac{\sin \theta
(s,t_{1})}{f(s,t_{1})}(\theta _{t}(s,t_{m}^{\prime })-\theta
_{t}(s,t_{1}))ds=0  \label{5.29}
\end{equation}%
\noindent by a similar reasoning as before. On the other hand, we estimate%
\begin{eqnarray*}
&&\frac{\sin \theta (s,t_{1})}{f(s,t_{1})}(\theta _{t}(s,t_{m}^{\prime
})-\theta _{t}(s,t_{1})) \\
&\geq &\frac{1}{2\max_{\bar{\Omega}_{1}}|f|}\frac{M^{\prime }}{2}\text{ (or }%
\leq \frac{-1}{2\max_{\bar{\Omega}_{1}}|f|}\frac{M^{\prime }}{2},\text{
resp.)}
\end{eqnarray*}%
\noindent by (\ref{5.27.3}), which contradicts (\ref{5.29}). Thus we have
shown that $\theta _{t}$ is continuous. Since $\theta _{s}$ is continuous by
Corollary C.1, we conclude that $\theta $ $\in $ $C^{1}$ in $s,$ $t$
coordinates, and hence $\theta $ $\in $ $C^{1}$ in $x,$ $y$ coordinates. It
follows that the characteristic curves are then $C^{2}$ smooth in view of (%
\ref{1.5}). Similarly the seed curves are also $C^{2}$ smooth.

To compute $\theta _{t},$ we first show that $D_{s}$ ($=$ $f^{-1}N^{\perp
}(D))$ exists and is continuous. Observe that 
\begin{equation}
x_{s}=\frac{\sin \theta }{f},x_{t}=\frac{\cos \theta }{gD},y_{s}=-\frac{\cos
\theta }{f},y_{t}=\frac{\sin \theta }{gD}  \label{5.34}
\end{equation}

\noindent and ($x_{s})_{t}$ (($y_{s})_{t}$ ,resp.) exists and is continuous
since $\theta $ $\in $ $C^{1}$ and $f$ is $C^{1}$ in $t$ by Theorem C. It
follows from Lemma 5.4 that ($x_{t})_{s}$ (($y_{t})_{s},$ resp.) exists and
equals ($x_{s})_{t}$ (($y_{s})_{t},$ resp.). So by (\ref{5.34}) $D_{s}$
exists and is continuous (either $\cos \theta $ $\neq $ $0$ or $\sin \theta $
$\neq $ $0)$. Now computing $x_{st}$ $=$ $x_{ts}$ and $y_{st}$ $=$ $y_{ts},$
we get%
\begin{equation}
\frac{(\cos \theta )\theta _{t}f-f_{t}\sin \theta }{f^{2}}=\frac{-(\sin
\theta )\theta _{s}gD-(gD)_{s}\cos \theta }{(gD)^{2}}  \label{5.35}
\end{equation}

\noindent and%
\begin{equation}
\frac{(\sin \theta )\theta _{t}f-f_{t}(-\cos \theta )}{f^{2}}=\frac{(\cos
\theta )\theta _{s}gD-(gD)_{s}\sin \theta }{(gD)^{2}}.  \label{5.36}
\end{equation}

\noindent Multiply (\ref{5.35}) by $\cos \theta $ and (\ref{5.36}) by $\sin
\theta ,$ respectively, and then add the resulting equalities to obtain%
\begin{equation}
\frac{\theta _{t}f}{f^{2}}=\frac{-(gD)_{s}}{(gD)^{2}}.  \label{5.37}
\end{equation}

\noindent We can now compute 
\begin{eqnarray*}
\theta _{t} &=&-\frac{fg_{s}}{g^{2}D}-\frac{fD_{s}}{gD^{2}} \\
&=&\frac{(rot\vec{F})g}{(gD)^{2}}-\frac{fD_{s}}{gD^{2}}\text{ (by (\ref{1.9}%
) and noting that }N^{\perp }g=fg_{s}) \\
&=&\frac{rot\vec{F}}{gD^{2}}-\frac{N^{\perp }(\log D)}{gD}
\end{eqnarray*}

\noindent which is (\ref{1.11}).

\endproof%

\bigskip

We remark that seed curves may only be $C^{1}$ smooth, but not $C^{2}$
smooth if $u$ is only Lipschitzian, but not $C^{1}$ smooth even in the case
of $H$ $=$ $0.$ An example is given by $u(x,y)$ $=$ $xy$ for $y$ $>$ $0$ and 
$u$ $=$ $0$ for $y$ $\leq $ $0$ with $\vec{F}$ $=$ $(-y,x)$ (see Example 7.2
in \cite{CHY})$.$ The seed curve $\Lambda _{a}$ passing through a point $(a,$
$0),$ $a$ $\neq $ $0,$ is the union of the straight line \{$x$ $=$ $\pm a\}$
for $y$ $\geq $ $0$ and the semi-circle of center $(0,$ $0)$ and radius 
\TEXTsymbol{\vert}$a|$ for $y$ $<$ $0.$ It is easy to see that $\Lambda _{a}$
is $C^{\infty }$ smooth except at $(\pm a,$ $0)$ where $\Lambda _{a}$ is
only $C^{1}$ smooth, but not $C^{2}$ smooth. Note that $u$ is not $C^{1}$
smooth at the $x$-axis while it is a Lipschitzian $p$-minimizer on any
bounded plane domain. For more examples, see (\cite{Ri}).

\end{document}